\documentclass[12pt,fleqn,a4paper]{amsart}

\usepackage[utf8]{inputenc}
\usepackage[in]{fullpage}
\usepackage[colorlinks,breaklinks,allcolors=blue]{hyperref}
\usepackage{amsfonts,amssymb,amsthm}
\usepackage[alphabetic,nobysame]{amsrefs}

\numberwithin{equation}{section}\swapnumbers
\usepackage{color}
\definecolor{changecolor}{rgb}{0,0.7,0}

\usepackage{parskip}

\usepackage[cmtip,matrix,arrow]{xy} \SelectTips{cm}{10}
\newcommand{\cxymatrix}[1]{\vcenter{\xymatrix@=15pt{#1}}}
\usepackage{rotating}

\newcommand{\xysubseteq}{\ar@{}[r]|{\displaystyle\subseteq}}
\newcommand{\xysubseteqdown}{\ar@{}[d]|{\rotatebox{90}{$\supseteq$}}}

\usepackage{tikz}
\usetikzlibrary{calc}

\usepackage{aliascnt}

\newtheorem{theorem}{Theorem}[section]
\newaliascnt{lemma}{theorem}

\newtheorem{lemma}[lemma]{Lemma}
\aliascntresetthe{lemma}
\newaliascnt{corollary}{theorem}

\newtheorem{corollary}[corollary]{Corollary}
\aliascntresetthe{corollary}
\newaliascnt{proposition}{theorem}

\newtheorem{conjecture}{Conjecture}

\newtheorem{proposition}[proposition]{Proposition}
\aliascntresetthe{proposition}

\theoremstyle{definition}
\newaliascnt{definition}{theorem}

\newtheorem{definition}[definition]{Definition}
\aliascntresetthe{definition}

\newaliascnt{remark}{theorem}
\newtheorem{remark}[remark]{Remark}
\aliascntresetthe{remark}
\newtheorem{remarks}[remark]{Remarks}

\newtheorem*{remark*}{Remark}

\newaliascnt{example}{theorem}

\aliascntresetthe{example}

\usepackage{enumitem}
\setlist[enumerate,2]{label=\textit{\alph*)},ref=\textit{\alph*})}
\setlist[enumerate,1]{label=\textit{\roman*)},ref=\textit{\roman*})}

\usepackage[capitalize]{cleveref}
\usepackage{microtype}

\newcommand{\Pq}{{\overline P}}
\renewcommand{\P}{{\mathbf P}}

\newcommand{\Xq}{{\overline X}}
\newcommand{\Yq}{{\overline Y}}
\newcommand{\Zq}{{\overline Z}}
\newcommand{\Gq}{{\overline G}}

\newcommand{\CC}{\mathbb{C}}

\newcommand{\QQ}{\mathbb{Q}}

\renewcommand{\rho}{\varrho}
\renewcommand{\phi}{\varphi}
\renewcommand{\epsilon}{\varepsilon}

\newcommand{\leer}{\varnothing}
 
\renewcommand{\>}{\rangle}
\newcommand{\into}{\hookrightarrow}

\renewcommand{\[}{\begin{equation}}
\renewcommand{\]}{\end{equation}}

\DeclareMathOperator{\rk}{rk}

\DeclareMathOperator{\Ad}{Ad}

\DeclareMathOperator{\md}{mod}
\DeclareMathOperator{\chr}{char}

\newcommand{\G}{{\mathbf{G}}}

\newcommand{\XS}{{\widetilde X}}
\newcommand{\YS}{{\widetilde Y}}

\newcommand{\yS}{{\widetilde y}}

\newcommand{\fB}{\mathfrak{B}}
\newcommand{\RR}{\mathbb{R}}
\newcommand{\Irr}{\mathfrak{Irr}}

\renewcommand{\mod}{/\!\!/}
\def\|#1|{\operatorname{#1}}

\newcommand{\ab}{\mathfrak{a}}
\newcommand{\ann}{\perp}

\newcommand{\codim}{\operatorname{codim}}

\newcommand{\g}{\mathfrak{g}}
\newcommand{\h}{\mathfrak{h}}

\newcommand{\kkk}{\Bbb K}
\newcommand{\lv}{\mathfrak{l}}

\newcommand{\st}{^{\rm st}}

\newcommand{\p}{\mathfrak{p}}
\newcommand{\m}{\mathfrak{m}}
\newcommand{\s}{\mathfrak{s}}

\newcommand{\tm}{\mathfrak{t}}

\newcommand{\gm}{\mathfrak{g}}

\newcommand{\Ru}[2][\relax]{{#2}_{{u}}\ifx#1\relax\else{\mspace{-7.8mu}}^{#1}\fi}

\newcommand{\K}{\mathcal{K}}

\newcommand{\dvr}{\operatorname{div}}

\newcommand{\T}[2][\relax]{T_{#2\ifx#1\relax\else,#1\fi}}
\newcommand{\CT}[2][\relax]{T^*_{#2\ifx#1\relax\else,#1\fi}}
\newcommand{\Ch}{\Lambda}

\newcommand{\Mom}[1][\relax]{\mu_{#1}}

\newcommand{\W}[1]{W_{#1}}
\newcommand{\WE}[1]{W_{(#1)}}

\newcommand{\CN}[2][\relax]{{\mathcal N}({#2\ifx#1\relax\else,#1\fi})}
\newcommand{\CNt}[2][\relax]{{\widetilde{\mathcal N}}({#2\ifx#1\relax\else,#1\fi})}

\newcommand{\CNB}[3][\relax]{N^*_{#2/#3\ifx#1\relax\else,#1\fi}}

\newcommand{\logd}[1][\relax]{\psi_{#1}}
\newcommand{\Logd}[1][\relax]{\Psi_{#1}}

\newcommand{\CTt}[1]{\widetilde{T}_{#1}}

\newcommand{\Gal}{\operatorname{Gal}}

\newcommand{\PU}{P_u}
\newcommand{\La}{L_{an}}
\newcommand{\Lv}{L}
\newcommand{\LXk}{L_{X,k}}
\newcommand{\AXk}{A_{X,k}}

\newcommand{\PXk}{P_{X,k}}

\newcommand{\LX}{L_{X}}

\newcommand{\ax}{\mathfrak a_{X}}
\newcommand{\PX}{P_{X}}
\newcommand{\PUX}{P_{X,u}}
\newcommand{\pux}{\mathfrak p_{X,u}}

\newcommand{\LXo}{L_{0}}
\newcommand{\AXo}{A_{0}}
\newcommand{\TXo}{T_{0}}
\newcommand{\axo}{\mathfrak a_{0}}
\newcommand{\LXS}{M_X}
\newcommand{\LS}{L(X)}
\newcommand{\ls}{{\mathfrak l(X)}}

\newcommand{\Id}{{\rm Id}}

\newcommand{\X}{X^\circ}

\newcommand{\el}{{\rm el}}
\newcommand{\an}{{\rm an}}

\newcommand{\Bal}{sPs}
\newcommand{\Ual}{ sP_us}

\newcommand{\mcd}{m_{cd}}   
\newcommand{\nd}{n_{cd}}  
\newcommand{\ld}{l_w} 
\newcommand{\eff}{{\rm eff}}

 \newcommand{\N}{\mathcal{N}}
 \newcommand{\NT}{\widetilde{\mathcal{N}}}

\newcommand{\muG}{\mu_X/\!\!/G}

\newcommand{\B}{\mathfrak{B}}     
\newcommand{\PS}{\mathfrak{B}_{0}}  

\newcommand{\ds}{{\raisebox{-.7\height}{$\cdot$}}\hskip -.65mm
{\raisebox{-.17\height}{$\cdot$}}\hskip -.65mm
{\raisebox{0.36\height}{$\cdot$}}\hskip -.65mm
{\raisebox{.9\height}{$\cdot$}}}

\DeclareMathOperator{\Rad}{\rm Rad}

\def\|#1|{\operatorname{#1}}

\newcommand{\DY}[1]{\Delta_{#1}}
\newcommand{\Chpr}[1]{\Ch({#1})^{\rm {pr}}}
\newcommand{\aY}[1]{\ab_{#1}}

\newcommand{\AF}[1]{A_{#1}}

\newcommand{\lvo}{\lv_{X,0}}
\newcommand{\lve}{\lv_{X,\rm{eff}}}
\newcommand{\LXoo}{L_{X,0}}

\newcommand{\iso}{ \stackrel{\sim}{\longrightarrow}}
\newcommand{\isol}{ \stackrel{\sim}{\longleftarrow}}

\title[]{Complexity of actions over perfect fields}

\author[]{Friedrich Knop}
\address[]{Dept. Mathematik, FAU Erlangen-Nürnberg, Cauerstraße 11,
D-91058 Erlangen}

\author[]{Vladimir S. Zhgoon} \address[]{Science Research Institute of
  System Studies, National Research University Higher School of
  Economics}

\begin{document}

\begin{abstract} Let $G$ be a connected reductive group over a perfect field $k$ acting on  an algebraic variety $X$ and let $P$ be a minimal parabolic subgroup of $G$.  
For $k$-spherical $G$-varieties  we prove finiteness result for $P$-orbits that contain $k$-points.  
This is a consequence of an equality on  $P$-complexities of $X$ and  of any $P$-invariant $k$-dense subvariety in $X$, 
which  generalizes a corresponding result of E.B.Vinberg in the case of algebraically closed field $k$. Also we introduce an action of the restricted Weyl group $W$ on the set of $k$-dense $P$-invariant closed  subvarieties of $X$
of maximal $P$-complexity and $k$-rank in the case of $\chr k =0$ and on the set of all $k$-dense $P$-orbits in the case of real spherical variety which generalizes the action on $B$-orbits introduced by F.Knop
in the algebraically closed field case.  We also introduce a  little Weyl group related with this action and describe its generators in terms of the generators  of $W$ which generalize the description of M.Brion in algebraically closed field case.  
\end{abstract}
\maketitle

\begin{flushright} \it {To the memory of  Èrnest Vinberg}
\end{flushright}

\section{Introduction}\label{s:Introduction}

Let $k$ be a ground field, let $G$ be a connected reductive algebraic
group defined over $k$, and let $X$ be a $G$-variety. For $k$
algebraically closed and of characteristic zero, the following theorem
was proved independently by M.Brion \cite{Brion}, and E.B.Vinberg
\cite{Vinberg} (the latter based on the paper of V.L.Popov
\cite{Popov}). Later T.Matsuki, \cite{Matsuki}, gave a completely
different and much simpler argument which also works in positive
characteristic.

\begin{theorem}\label{thm:VinbergBrion}
  Assume that a Borel subgroup $B$ of $G$ has an open orbit in $X$ (in
  which case $X$ is called spherical). Then the number of $B$-orbits
  in $X$ is finite.
\end{theorem}

The theorem is sharp in two ways. First, it does not hold if $B$ is
replaced by $G$ or any other non-solvable parabolic subgroup of
$G$. Secondly, the theorem is not true for $B$-actions which are not
the restriction of a $G$-action.

The aim of this paper is to generalize the \cref{thm:VinbergBrion}
above to more general ground fields. The obvious way is to replace the
Borel subgroup which may not be defined over $k$ by a minimal
parabolic $k$-subgroup $P$. However this is not enough. Masuki
\cite{Matsuki}*{Rmk.~7} gave for $k=\RR$ an example of a $k$-variety
$X$ having infinitely many $P$-orbits despite having an open
$P$-orbit. Most of these orbits do not contain $k$-points, though. So,
Matsuki conjectured \cite{Matsuki}*{Conj.~2} that in the presence of
an open $P$-orbit, $X(\RR)$ has only finitely many $P(\RR)$-orbits. In
the same paper, Matsuki showed how to reduce his conjecture to the
case of groups of real rank one. For that case, Kimelfeld
\cite{Kimelfeld} had already classified all relevant spaces three
years earlier but his paper was overlooked. Matsuki's conjecture was
subsequently proved independently by Bien \cite{Bien} and
Krötz-Schlichtkrull \cite{KSch}.

It is now tempting to replace in Matsuki's conjecture $\RR$ by an
arbitrary field $k$ but in that form it is not true. The reason is
simply that there may be $P$-orbits $Y$ where $Y(k)$ decomposes into
infinitely many $P(k)$-orbits. The easiest example is $k=\QQ$,
$G=X=\G_m$ with ``square'' action $a*x:=a^2x$. Here
$X(k)/P(k)=\QQ^\times/(\QQ^\times)^2$ is infinite. This difficulty
is overcome by only considering $P$-orbits which contain a
$k$-point. Now our generalization of \cref{thm:VinbergBrion} is:

\begin{theorem}\label{thm:MainSpherical}
  Let $k$ be an infinite perfect field and assume that a minimal
  parabolic $k$-subgroup $P$ of $G$ has an open orbit in $X$ (in which
  case $X$ is called $k$-spherical). Then the number of $P$-orbits $Y$
  of $X$ with $Y(k)\ne\emptyset$ is finite.
\end{theorem}

The theorem also obviously also true for finite fields but in that
case \cref{thm:VinbergBrion} is stronger (since then $P=B$). So we
excluded that case. The perfectness assumption comes from Kempf's
instability theorem \cite{Kempf}. We don't know whether the theorem
holds for non-perfect fields.

Borel-Serre \cite{BorelSerre} have proved that over local fields of
characteristic zero for every $P$-orbit $Y$ the set $Y(k)/P(k)$ is
finite. This immediately implies Matsuki's conjecture for local fields:

\begin{corollary}
  Let $k$ be a local field of characteristic zero and let $X$ be a
  $k$-spherical $G$-variety. Then $X(k)/P(k)$ is finite.
\end{corollary}

Vinberg \cite{Vinberg} has proved an even more general
version of \cref{thm:VinbergBrion} which works for any $G$-variety. To
this end, he introduced the complexity $c(Y)$ of a $B$-invariant
subvariety $Y$ as the dimension of a birational quotient $Y/B$. So
$c(Y)=0$ means that $B$ has an orbit in $Y$. For this notion, Vinberg
proved $c(Y)\leq c(X)$ for any $B$-stable subvariety.
\cref{thm:VinbergBrion} is just the special case of $c(X)=0$. In this
paper we prove an analogous extension of \cref{thm:MainSpherical}.
The idea is to consider $P$-stable subvarieties $Y$ for which
$Y(k)$ is dense. For the precise statement see \cref{MAIN_COMP} below.

The original proofs of Vinberg and Brion used a reduction argument
where $X$ is deformed to a horospherical variety. We do not know how
to adapt this to our situation since it seems to be impossible to
control $k$-points under a deformation of $X$. Instead we use the
aforementioned reduction argument of Matsuki to semisimple
groups. Over algebraically closed fields there is basically only one,
namely $SL(2,k)$. In general, the structure of semisimple
rank-1-groups is much more involved so that we cannot simply resort to
a classification. Instead, we use some general structure theorems for
the action of anisotropic groups which may interesting in their own
right.

After proving the finiteness theorem it suggests itself to generalize
also the results of \cite{BORB}. Assume for simplicity that $X$ is
$k$-spherical and let $\fB(X)$ be the finite set of $P$-orbits
$Y\subseteq X$ such that $Y(k)\ne\emptyset$. In the setting of $k$
being algebraically closed fields of characteristic not $2$, a
canonical action of the Weyl group $W$ of $G$ on $\fB(X)$ has been
constructed. For general fields, we replace $W$ by the restricted Weyl
group $W_k$. For every simple reflection of $W_k$ be define an
involution of $\fB(X)$ and we conjecture that these involutions extend
to an action of $W_k$ on $\fB(X)$. The problem is of course the
verification of the braid relations.

Here we have only partial results. Apart from algebraically closed
fields we show that the conjecture holds for the field $k=\RR$ of
real numbers. The proof uses the same Hecke algebra approach as in
\cite{BORB}.

Otherwise, we show that for ${\rm char}\,k=0$, the $W_k$-action exists
on a certain subset $\fB_0(X)$ of $\fB(X)$ where $\fB_0(X)$ consists
of the $P$-orbits of maximal rank or, equivalently, of those
$P$-orbits $Y$ for which $\dim Y/\Rad_uP$ is maximal. This subset
contains the open $P$-orbit $(X)$ and we show that
$W_k\circ (X)=\fB_0(X)$. Moreover, we prove that the stabilizer
$W_{(X)}$ of the open $P$-orbit is related to the little Weyl group
$W_k(X)$ of $X$ as defined in \cite{KK}. To this end we use, as in
\cite{BORB}, the momentum map on the cotangent bundle of $X$.

Finally, we also generalize Brion's approach \cite{BBORB} for
describing generators of $W_{(X)}$. These are either reflections or
products of commuting reflections in $W_k$ .


We indicate the structure of the paper. In \cref{s:rank-1-case} we
study the case where the semisimple $k$-rank of $G$ is
one. 
\cref{s:generalcase} is dedicated to
deriving an inequality for complexities in the general case. In
\cref{s:realcase} we consider Hecke algebra approach to the action
of the Weyl group on the set of $P$-orbits of a real spherical
variety. In \cref{s:principal} we study equivariant geometry of
cotangent vector bundle and it relation with the action of a Weyl
group on the $k$-dense sheets of $P$-orbits with maximal complexity,
rank and homogeneity. \cref{s:functoriality} is dedicated to comparison of
the actions of the Weyl groups under field extensions.
\cref{s:generators} is dedicated to the proof of the fact that
little Weyl group is generated by the set of reflections in the roots
of $G$ and by the products of commuting reflections. In the
\cref{s:monodromy} for a real homogeneous spherical variety that have
wonderful compactification we provide the action of fundamental group
of an associated flag variety (a unique closed orbit in the wonderful
compactification) on the set of connected components of the open
$P$-orbit (which are $P(\RR)^\circ$-orbits). This action is
important since it preserves the set of $G(\RR)$-orbits.

 One of the natural questions that arise in relation with our work is studying  the possibility of extention of  our action
 of the restricted Weyl group in the real case on the set of principal orbits to the action on  to the set of  $B(\Bbb R)$-orbits in $X(\Bbb R)$.
  S.Cupit-Foutou and D.A.Timashev  first raised this question  in \cite{TC} in the  case of
  spherical varieties of the real split groups  and provided the extension of a subgroup of the Weyl group under some technical assumptions on the class of varieties.
  These assumptions  were necessary, since  they also discovered an example where the generators of their action generate a group which  is different from the
  Weyl group.

Finally, let us mention that our work has started during Vinberg's
$80^{\rm th}$ anniversary conference in Moscow 2017 and was announced
officially on Brion's $60^{\rm th}$ anniversary conference in Lyon
2018. We would like express our warmest wishes and greatest respect to Michel Brion.
Unfortunately,  Èrnest  Vinberg had passed away on 12th May 2020, so we wish to dedicate this paper to his blessed memory.

\section{Complexities}\label{s:Complexities}

Let $k$ be an infinite perfect field with algebraic closure $K$. A
(not necessarily irreducible) $k$-variety $X$ is called
\emph{$k$-dense} if $X(k)$ is Zariski dense in $X$. Note that for
$k$-dense varieties, irreducibility is equivalent to absolute
irreducibility.

Let $H$ be a linear algebraic $k$-group. An $H$-variety $X$ is called
\emph{$H$-irreducible} 
if $H$-span of the set of $k$-points is Zariski dense in $X$, i.e. $\overline{HX(k)}=X$.
 When $H$ is connected then
$H$-irreducibility is equivalent to irreducibility.

Let $k(X)$ be the total ring of fractions of $X$. If $X$ is
$H$-irreducible then $k(X)^H$ is a field and we define the
\emph{$H$-complexity of $X$} as
\[
  c(X/H):=\operatorname{trdeg}_kk(X)^H.
\]
The notation ``$X/H$'' is purely formal and does not mean that an
orbit space $X/H$ exists. If it does, though, then $c(X/H)=\dim X/H$.
If $X^\circ\subseteq X$ is any irreducible component and
$H^\circ\subseteq H$ the component of unity then clearly
$c(X/H)=c(X^\circ/H^\circ)$.

The most important instance of $H$-complexity is as follows: Assume
$X$ is a $k$-dense, irreducible $G$-variety where $G$ is a connected
reductive $k$-group. Let $P\subseteq G$ be a be a minimal parabolic
$k$-subgroup. Then the \emph{$k$-complexity of $X$} is defined as its
$P$-complexity:
\[
  c_k(X):=c(X/P).
\]

By $\fB(X)$ let us denote the set of $k$-dense closed subvarieties $Y\subset X$ such that $c_k(Y)=c_k(X)$.
If $X$ is $k$-spherical we shall prove that $\fB(X)$ coincides with the set of $P$-orbit
closures $\overline{Px}$ where $x\in X(k)$.

More generally, given an $H$-variety $X$, we are interested in the set
$\Irr(X/H)$ of all closed $H$-irreducible subvarieties
$Z\subseteq X$.

\begin{lemma}\label{P_G_H}

  Let $G$ be a connected reductive group, let $P\subseteq G$ be a
 $k$-parabolic subgroup, and $H\subseteq G$ a closed $k$-subgroup. Put
  $X:=G/H$ and $Y:=G/P$ inducing the diagram
  
  \[
    \cxymatrix{
      &G\ar[dl]_{\pi_H}\ar[dr]^{\pi_P}\\
      G/H&&G/P }
  \]
  with $\pi_H(g)=gH$ and $\pi_P(g)=g^{-1}P$. Then there is an
  injective inclusion preserving map
  \[
    \Irr(Y/H)\to\Irr(X/P):Z\mapsto Z':=\pi_H(\pi_P^{-1}(Z))
  \]
  which also preserves complexities:
  \[
c(Z/H)=c(Z'/P).
    \]
\end{lemma}

\begin{proof}

  The problem is that the preimage of a $k$-dense subvariety is, in
  general, not $k$-dense. But this is true for $\pi_P$ since
  $G(k)\to Y(k)$ is surjective and $P$ is $k$-dense. This implies that
  $\pi_P^{-1}$ induces a bijection
  $\Irr(Y/H)\overset\sim\to\Irr(G/(H\times
  P))$. The latter maps via $\pi_H$ into $\Irr(X/P)$. Since this map preserves codimension and relative  codimension of the subsets, it is
  clear that complexities are preserved.
\end{proof}

\begin{remarks} 1. Observe that the case $c=0$ means that the
  correspondence maps orbit closures to orbit closures.

  2. Observe that $\Irr(Y/H)$ and its inclusion into
  $\Irr(X/P)$ depends only on $G(k)$-conjugacy class of $H$. It
  does not just depend on $X$, though. A typical example would be
  $k=\mathbb{R}$ and $X=SL(3)/SO(3)$, the space of unimodular
  quadratic forms. Then $X(\mathbb{R})$ has two $G(\mathbb{R})$-orbits
  with isotropy groups $H_1=SO(3,\RR)$ and $H_2=SO(1,2,\RR)$,
  respectively. Since $H_1$ acts transitively on $Y(\RR)$ we have
  $\Irr(Y/H_1)=\{Y\}$. On the other hand, $H_2(\CC)$ has four
  orbits in $Y(\CC)$ and all of them are $\RR$-dense. Thus, $|\Irr(Y/H_2)|=4$.

  \end{remarks}
  
\section{The rank-1-case}\label{s:rank-1-case}

In the following, $k$ is an arbitrary infinite perfect field. We start
with a codimension-$1$-variant of a well-known theorem of Rosenlicht.

\begin{lemma}\label{lemma:pseudoquotient}

  In the following, everything should be defined over $k$. Let $P$ be
  a connected group and let $R\subseteq P$ a connected normal
  subgroup. Let $X$ be a smooth $P$-variety such that all stabilizers
  $R_x$ are finite. Let, moreover, $Y\subset X$ be a $P$-stable prime
  divisor. Then there is a $P$-stable open subset $X_0\subseteq X$
  with $X_0\cap Y\ne\leer$, a smooth $P$-variety $Z$, and a surjective
  $R$-invariant $P$-morphism $\phi:X_0\to Z$ such that every fiber is
  a finite union of $R$-orbits.

\end{lemma}

\begin{proof}

  By Rosenlicht's theorem \cite{Rosenlicht}*{Thm.~2}, there is a
  $P$-stable open subset $X_1\subseteq X$ such that the orbit space
  $Z_1=X_1/R$ exists. We may assume that $Z_1$ is smooth. If
  $X_1\cap Y\ne\leer$ we are done. So assume that $X_1\cap Y=\leer$.
  The divisor $Y$ induces a $P$-invariant valuation $v_Y$ of
  $K(X)$. Its restriction $v'=v_Y|_{K(Z_1)}$ is $P/R$-invariant and
  induces an equivariant open embedding $Z_1\subseteq Z=Z_1\cup Y'$
  where $Z$ is a smooth $P$-variety and $Y'\subset Z$ is a
  $P$-invariant prime divisor. The morphism $X_1\to Z_1$ extends, by
  construction, to a $P$-equivariant $k$-morphism $\phi:X_0\to Z$
  where $X_0\subseteq X$ is open $P$-stable with $X_0\supseteq X_1$
  and $Y_0:=X_0\cap Y\ne\leer$. Since the map $Y_0\to Y'$ between
  divisors is dominant (again by construction), we may shrink $Z$ such
  that $\phi$ is surjective and equidimensional. By assumption, all
  $R$-orbits in $X$ have dimension $\dim R$ which is also the
  dimension of the fibers of $\phi$. The $R$-invariance of $\phi$
  implies that every fiber is a finite union of $R$-orbits.
\end{proof}

Next we extend a well-known property of torus actions.

\begin{lemma}\label{lemma:affineNbhd}

  Let $G$ be an elementary $k$-group acting on a locally linear
  $k$-variety $X$. Then every $x\in X(k)$ is contained in a
  $G$-stable, affine, open subset $X_0\subseteq X$.

\end{lemma}

\begin{proof}

  We have $G=MA$ with $M$ anisotropic and $A$ a split central
  torus. By local linearity, we may assume that $X$ is a subvariety of
  $\P(V)$ where $V$ is a $G$-module. Let $Z:=\Xq\setminus X$ be the
  boundary and let $\tilde Z\subseteq V$ be the affine cone over
  $Z$. Let, moreover, $0\ne\tilde x\in V(k)$ be a lift of $x$. Since
  $M$ is anisotropic, the orbit $M\tilde x$ is closed in $V$ by
  Kempf's theorem (see \cite{Kempf}*{Remark after Cor.\ 4.4}) and does
  not meet $\tilde Z$. It follows that there is a homogeneous
  $M$-invariant $f$ on $V$ which vanishes on $\tilde Z$ and with
  $f(M\tilde x)=1$. Let $f=\sum_\chi f_\chi$ be the decomposition of
  $f$ into $A$-eigenfunctions. Since $Z$ is $G$-stable, all $f_\chi$
  are $G$-eigenfunction which vanish on $\tilde Z$. One of them, say
  $f_{\chi_0}$, is non-zero in $\tilde x$. Then
  $X_0=X\setminus\{f_{\chi_0}=0\}$ is clearly a $G$-stable, affine,
  open subset of $X$ containing $x$.
\end{proof}

Let $P$ be a quasi-elementary group with Iwasawa decomposition $P=\La AP_u$, where $\La$ is maximal anisotropic subgroup, $P_u$ is the unipotent radical and $A$ is the maximal
split torus. 
For a $P$-invariant subset $Y$. By shrinking the rational quotient $Y/AP_u$ we can assume that it is smooth 
and by previous theorem we may assume that  is affine. In the case of $\chr k=0$ by Richardson theorem (see \cite{inv}) there exists the  stabilizer of general position for $\La$ on  
 rational quotient $Y/AP_u$  that we will denote by $M_Y$. 
 By $s_k(Y)$ let us denote the dimension of $\La$-orbit on $Y/AP_u$ of general position. 

\begin{lemma}\label{lemma:anisotropic}

  Let $P$ be a quasi-elementary group (e.g. a minimal parabolic)
  acting on a smooth $k$-dense variety $X$. Let $Y\subset X$ be a
  proper $P$-stable subvariety such that the set of $k$-points
  $x\in Y(k)$ with anisotropic stabilizer $P_x$ is dense in $Y$. Then
  $c(Y/P)<c(X/P)$,  $\rk_k(X)\ge \rk_k(Px)$ and $s_k(Y)\le s_k(X)$.

\end{lemma}

\begin{proof}
  
  By removing $Y_{\|sing|}$ from $X$ and $Y$, we may assume that $Y$
  is smooth, as well.

  Let $\pi:\XS\to X$ be the blow-up of $X$ in $Y$ and let
  $\YS=\pi^{-1}(Y)$ be the exceptional divisor. Then $\XS$ and $\YS$
  are smooth and every fiber of $\pi|_D: \YS\to Y$ over a $k$-point is
  a projective space. If $\yS\in\YS$ we have $P_{\yS}\subseteq
  P_y$. This implies that also $\YS$ has a dense set of $k$-points
  with an anisotropic stabilizer. Since $c(\YS/P)\ge c(Y/P)$ it
  suffices to prove  first assertion for $(\XS,\YS)$. 
  Since we have a rational map between $\YS/AP_u\rightarrow Y/AP_u$ we also have $s_k(\YS)\ge s_k(Y)$.
  This reduces the problem to the case when $Y$ is of codimension $1$ in $X$.

  Let $R\subseteq P$ be the split radical of $P$, i.e., $R$ is normal,
  split solvable and $M=P/R$ is anisotropic. Since $R$ does not
  contain any nontrivial anisotropic subgroup, we have that $R_x$ is
  finite for every $x\in Y(k)$ with $P_x$ anisotropic. These points
  form a dense subset of $Y$. It follows that $R$ acts locally freely
  (i.e. with finite stabilizers) on an open subset $X'$ of $X$ with
  $X'\cap Y\ne\leer$. Replacing $X$ by $X'$ we may assume that the
  action of $R$ on $X$ is locally free.

  Now let $X_0\to Z$ be as \cref{lemma:pseudoquotient} and let
  $Y'\subseteq Z$ be the image of $X_0\cap Y$ which is an $M$-stable
  prime divisor. Since clearly $c(X/P)=c(Z/M)$ and
  $c(Y/P)=c(Y'/M)$. So after replacing $(X,Y)$ with $(Z,Y')$ we may
  assume that $P=M$ is anisotropic.

  Now \cref{lemma:affineNbhd} allows us to replace $X$ by an affine
  open subset. Since $Mx$ is closed in $X$ for $x\in X(k)$ we conclude
  that the generic orbits in $X$ and $Y$ are closed. Hence, $Y\mod M$
  is a proper subvariety of $X\mod M$ which implies
  $c(Y/P)=\dim Y\mod M<\dim X\mod M=c(X/P)$.
 The inequality  $s_k(X)\ge s_k(\YS)$  now follows from the fact that $\YS/P_uA$ is a subset of $X/P_uA$ and from  the semicontinuity of the dimension of $\La$-orbits.
 The same argument applied to $A$-orbits on the quotient by $P_u\La$ gives $\rk_k(X)\ge\rk_k(Px)$. 
\end{proof}

 \begin{corollary}\label{dense_orb_anis}
   Let $P$ be a quasi-elementary $k$-group acting on a smooth $k$-variety $X$. Assume
   that $X$ contains an open $P$-orbit $X_0$. Then the stabilizer
   $P_x$ of every point $x\in X(k)\setminus X_0(k)$ is isotropic.
 \end{corollary}

In the rest of this section, let $H^\circ$ be the \emph{reduced},
connected component of an algebraic group $H$.

\begin{lemma}\label{lemma:isolated}

  Let $G$ be a reductive group and $P=\La A\PU\subseteq G$ a parabolic
  subgroup. Let $X$ be a homogeneous $G$-variety. Let $Z\subseteq X$
  be the union of all closed $P$-orbits. Then:

  \begin{enumerate}

  \item\label{it:1.1} $x\in Z$ if and only if $P_x^\circ$ is parabolic
    in $G_x^\circ$.

  \item\label{it:1.2} $Z$ is closed in $X$.

  \item\label{it:1.3} Let $x\in Z$ and let $T_0\subseteq P_x$ be a
    maximal torus. Assume ${\PU}^{T_0}=\{1\}$. Then $Px$ is a connected
    component of $Z$.

  \end{enumerate}

\end{lemma}

\begin{proof}

  Let $H:=G_x$. Then $X=G/H$. Now consider the diagram
  \[
    \cxymatrix{
      &G\ar[dl]_{\pi_H}\ar[dr]^{\pi_P}\\
      G/H&&G/P }
  \]
  where $\pi_H(g)=gH$ and $\pi_P(g)=g^{-1}P$. Then
  $Z\mapsto\pi_H\pi_P^{-1}(Z)$ furnished an inclusion preserving
  bijection between closed $H$-stable subsets of $G/P$ and closed
  $P$-stable subsets of $G/H$. In particular, closed $H$-orbits
  correspond to closed $P$-orbits. Moreover, $Px\subseteq X$
  corresponds to $Hy\subseteq G/P$ where $y=eP$.

  \ref{it:1.1} Since $G/P$ is complete, $Hy$ is closed if and only if
  it is complete if and only if $H_y^\circ=(H\cap P)^\circ=P_x^\circ$
  is parabolic in $H^\circ=G_x^\circ$.

  \ref{it:1.2} Let $B_H\subseteq H^\circ$ be a Borel subgroup. Then an
  $H$-orbit is closed if and only if it contains a $B_H$-fixed
  point. Thus, the union of closed $H$-orbits is $H(G/P)^{B_H}$ which
  is closed since it is the image of the proper map
  $H/B_H\times(G/P)^{B_H}\to G/P$.

  \ref{it:1.3} The torus $T_0$ is contained in a Levi complement $L$
  of $P$. Let $P^-$ be opposite to $P$ with respect to $L$ and let
  $\PU^-$ be its unipotent radical. Then, as an $L$-variety, an open
  neighborhood of $y$ in $G/P$ is isomorphic to $\PU^-$. Since
  ${\PU}^{T_0}=1$ also $({\PU}^-)^{T_0}=1$. This means that $y$ is an isolated
  $T_0$-fixed point. Since $T_0\subseteq P_x=H\cap P$ we may choose
  $B_H$ such that $T_0\subseteq B_H$. Then, a fortiori, $y$ is an
  isolated point of $(G/P)^{B_H}$, as well. Since each closed
  $H$-orbit being a flag variety contains a unique $B_H$-fixed point
  the $H$-spans of different irreducible components of $(G/P)^{B_H}$
  are closed and do not intersect.  Therefore $Hy$ and $Px$ are open
  and closed in their respective union of all closed orbits.
\end{proof}

A $P$-orbit as in \ref{it:1.3} will be called
\emph{isolated}. Observe, that there exist only finitely many of them.

\begin{lemma}\label{lemma:affine}

  Let $G$ be a semisimple group with $\rk_kG=1$ and with minimal
  parabolic $P\subseteq G$. Let $X$ be an affine homogeneous
  $G$-variety and let $x\in X(k)$ a be $k$-point whose stabilizer
  $P_x$ is isotropic. Then $Px$ is isolated and $\rk_k(X)>\rk_k(Px)$.

\end{lemma}

\begin{proof}

  Since $X$ is affine, the stabilizer $H:=G_x$ is reductive, as
  well. It contains the isotropic subgroup $P_x$ so is also isotropic.
  Therefore $\rk_kH=1$ and $H$ contains a split torus $A_0$ of
  dimension $1$. Let $\lambda:\G_m\to A_0$ be an isomorphism. Then
  $\lambda$ induces parabolic subgroups
  $P_H^\pm:=P_H^\pm(\lambda)\subseteq H^\circ$ and
  $P_G^\pm:=P_G^\pm(\lambda)\subseteq G$ with
  $P_H^\pm=P_G^\pm\cap H^\circ$. Observe that $P_H^\pm$ and $P_G^\pm$
  are minimal parabolics of $H$ and $G$, respectively.

  We claim that the stabilizer $P_x^\circ=(H\cap P)^\circ$ is
  parabolic in $H$.  Since it is isotropic by assumption, it
  contains a subgroup $S$ which is either isomorphic to $\G_a$ or
  $\G_m$. By replacing $A_0$ by $hA_0 h^{-1}$ with a suitable
  $h\in H(k)$ we may assume that $S$ either lies in $R_u(P_H)$ or is
  equal to $A_0$.

  Let $Y:=G/P_G$ and $y_0:=eP_G\in Y(k)$. From the Bruhat
  decomposition
  \[\label{eq:Bruhat}
    Y(k)=\{y_0\}\cup P_G(k)sy_0
  \]
  we see that the set of $S$-fixed points in $Y(k)$ is either
  $\{y_0\}$ or $\{y_0,sy_0\}$. On the other hand, there is $g\in G(k)$
  with $P=gP_Gg^{-1}$. Then $P_x=H\cap P$ is equal to $H_y$ with
  $y:=gP_G\in Y(k)$. From $S\subseteq H_y$ we get
  $y\in\{y_0,sy_0\}$. But then $H^\circ y\cong H^\circ/P_H^\pm$ is
  closed in $Y$ which shows that $P_x^\circ=H_y^\circ=P_H^\pm$ is
  parabolic in $H$.

  Finally, since $P_x^\circ$ is parabolic in $H^\circ$ it contains a
  split torus $A\cong\G_m$ which acts with non-zero weights on
  $\PU$. The same is true for a maximal torus $T_0$ of $P_x$. This shows
  ${\PU}^{T_0}=1$ and we conclude with \cref{lemma:isolated}.
The proof of the inequality for the rank is postponed to the section where we calculate the rank lattices.
%
\end{proof}

Following Serre, a subgroup $H$ of $G$ is called
\emph{irreducible} if is not contained in a proper parabolic subgroup
of $G$.  In our settings we shall use more general definition  and reserve a different name for it:

\begin{definition}

  A $G$-variety over $k$ is \emph{homogeneously irreducible over k} if there is no
  $G$-morphism $\phi:X_0\to Z$ with $X_0\subseteq X$ open $G$-stable
  and $Z$ complete, homogeneous, and positive dimensional where
  $\phi$, $X_0$, and $Z$ are all defined over $k$. Otherwise, $X$ will
  be called homogeneously reducible.
\end{definition}

\begin{lemma}\label{lemma:maptoG/P}

  Let $G$ be a connected reductive $k$-group acting on a $k$-dense
  variety $X$. Assume that $X$ is homogeneously irreducible over $k$. Then the generic
  orbits of $X$ are affine.

\end{lemma}

\begin{proof}

  By Sumihiro's theorem, we may assume that $X$ is a subvariety of
  $\P(V)$ where $V$ is a finite dimensional $G$-module. Clearly it
  suffices to prove the assertion for the closure of $X$. Thus, we may
  assume that $X$ is closed. Let $\hat X\subseteq V$ be the affine
  cone of $X$.

  Suppose, for a generic point $x\in X$ the orbit $Gx$ is not affine
  or, equivalently, that the stabilizer $G_x$ is not reductive. Let
  $G'$ be the derived subgroup of $G$. Then also $G'_x$ is not
  reductive. Let $\hat x\in \hat X$ be a lift of $x$. Then the
  stabilizer $G'_{\hat x}$ can't be reductive, either. This in turn
  implies that the orbit $G'\hat x$ is not closed in $\hat X$. Thus, a
  generic $G'$-orbit in $\hat X$ is not closed.

  Let $m$ be the dimension of a generic $G'$-orbit in $\hat X$. Then
  the union $Y$ of all $G'$-orbits of dimension $<m$ is closed. It is
  a proper subset of $\hat X$ and has, by the above, the property that
  $Y$ meets the closure of a generic $G'$-orbit.

  Let $W'\subseteq k[V]$ be a finite dimensional rational submodule
  generating the ideal $\mathcal I(Y)$. Observe that, since $Y$ is a
  cone, $W'$ can be chosen to be homogeneous. Then the inclusion
  $W'\into k[V]$ induces a homogeneous $G$-morphism $\pi:V\to W$ where
  $W:=(W')^*$ with the property that $\pi^{-1}(0)=Y$.

  Let $X'\subseteq W$ be the closure of $\pi(\hat X)$. Then, by
  construction, the generic $G'$-orbit of $X'$ contains $0$ in its
  closure. Thus, Kempf's instability theory applies: there is a
  $1$-parameter subgroup $\lambda:\G_m\to G'$ which is optimal for the
  generic $G'$-orbit in $X'$. More precisely, let $X^+$ be the set of
  points $x\in X'$ such that $\lim_{t\to0}\lambda(t)x=0$. Then
  $G'X^+=X$. Moreover, $X^+$ is stable under $P:=P(\lambda)$, the
  parabolic of $G$ attached to $\lambda$ inducing a surjective
  morphism $\psi:G*_PX^+\to X'$  \footnote{where for the groups $G\supset H$ and a quasiprojective variety $Z$ by $G*_HZ$ we denote the quotient of $G\times Z$ by the $H$-action defined by $(g,z)\rightarrow (gh^{-1},hz)$.
  There is a  projection  $G*_HZ\rightarrow G/H$ given by $(g,z)\rightarrow gH$. If $Z$ is $H$-invariant subset of $G$-variety $X$ there is an action map $\psi:G*_HZ\rightarrow X$
  given by $\psi(g,z)=gz$ whose image is $GZ$.} Observe that $P\ne G$ since
  $\lambda(\G_m)$ must be a non-trivial subgroup of $G'$.

  The main property, though, is that $G'_x\subseteq P\cap G'$ for all
  generic $x\in X^+$. This immediately implies that also
  $G_x\subseteq P$ which means that $\psi$ is generically bijective,
  hence purely inseparable. Let $G_n\subseteq G$ be the $n$-th
  Frobenius kernel and $Z:=G/G_nP$. Then $Z$ is a complete $G$-variety
  such that for $n>\!\!>0$ there is a rational morphism $\pi'$
  completing the following square:
  \[
    \cxymatrix{G*_{P(\lambda)}X^+\ar[r]^>>>\psi\ar[d]&X'\ar@{.>}[d]^{\pi'}\\
      G/P\ar[r]&Z}
  \]
  Observe that $\pi'$ is defined over $k$ since $X'$ is $k$-dense and
  therefore $\lambda$ is Galois invariant.

  Finally, composition with $\hat X\to X'$ yields a rational map
  $\xymatrix{\hat X\ar@{.>}[r]&Z}$. This morphism is invariant for the
  scalar $\G_m$-action on $\hat X$. Therefore it descends to a
  rational $G$-morphism $\xymatrix{X\ar@{.>}[r]&Z}$.
\end{proof}

\begin{proposition}\label{prop:main}

  Let $G$ be a connected reductive group with $\rk_k^{\rm ss}G=1$ and
  with minimal parabolic subgroup $P$. Let $X$ be a $G$-variety and
  let $Y\subset X$ be a $k$-dense $P$-stable proper subvariety such
  that $GY$ is dense in $X$. Then 
   \begin{enumerate}

  \item\label{it: 1}
   either $c(Y/P)<c_k(X)$ or $c(Y/P)=c_k(X)=c(X/G)$. 
   \item\label{it: 2}
  $\rk_k(Y)\leq \rk_k(X)$, if moreover $\rk_k(Y)=\rk_k(X)$, then $X$ is $k$-horospherical.
   \item\label{it:3}
  The inequality $s_k(Y)>s_k(X)$ can hold only in the case when $X$ is homogeneously irreducible and in this case
  $\rk_k(Y)<\rk_k(X)$. 
    \end{enumerate}
\end{proposition}

\begin{corollary} In the lexicographic order we have $$(c_k(Y),\rk_k(Y), s_k(Y))\leq (c_k(X),\rk_k(X), s_k(X)).$$
\end{corollary}

\begin{remark} The inequality $s_k(Y)>s_k(X)$ can occur. This   happens   for suitable $n,m$ when $Y$ is the closed $P$-orbit in a homogeneous space $X=SO(1,n)/SO(1,m)\times SO(n-m)$ over $k=\Bbb R$.
However in this case $\rk_k(Y)<\rk_k(X)$.
\end{remark}

\begin{definition} The closed $P$-variety $Y$ (the family of generic $P$-orbits in $Y$) is called principal if $(c_k(Y),\rk_k(Y), s_k(Y))$ is maximal in the lexicographic order. We denote the set of principal varieties 
by $\PS(X)$.
\end{definition}


\begin{proof}

  Let $Z(G)\subseteq G$ be the center of $G$. By Rosenlicht, $X$
  contains an open $G$-stable subset $X^0$ such that the orbit spaces
  $X^0/G$ and $X^0/Z(G)$ exist. It follows from $X=GY$ that
  $Y^0=X^0\cap Y\ne\leer$ and that $Y^0/G$ exists. Moreover, dividing
  by $Z(G)$ does not change complexities. So we may replace $G$, $P$
  $X$, $Y$ by $G/Z(G)$, $P/Z(G)$, $X^0$, $Y^0$ (or any smaller open
  subset) and assume that $G$ is semisimple of $k$-rank $1$ and that
  $X/G$ exists.


  Since $X=GY$, the map $Y\to X/G$ is surjective.  Assume $c(Y/P)\ge  c_k(X) \ge c(X/G)$. 
  
  Assume first that the generic orbit of $X$ is affine. After
  shrinking $X$, if necessary, we may assume that all orbits are
  affine and that the union $Z\subseteq X$ of all closed $P$-orbits is
  closed. Let $Z_0\subseteq Z$ be the union of all components such
  that $Z_0/P\to X/G$.  After shrinking $X$, we may assume that the intersection
  $Z_0\cap Gx$  with every orbit is the union of all isolated $P$-orbits of $Gx$.  
  Hence all $x\in X(k)$ with isotropic stabilizer
  $P_x$ lie in $Z_0$ by \cref{lemma:affine}. 
  If the set of points  $x\in Y(k)$ with an anisotropic stabilizer $P_x$ is dense in
  $Y$ by \cref{lemma:anisotropic} we get $c(Y/P)< c_k(X)$.
  Otherwise $Y$  is the closure of the set of points  $x\in Y(k)$ with isotropic stabilizer and $Y\subset Z_0$. Then
  $c(Y/P)\le c(Z_0/P)=c(X/G) \le c(X/P)$ and the equality can hold iff $P$ has an open k-dense orbit in $Gx$.
  By Corollary \ref{dense_orb_anis} there is only one $k$-dense orbit with an anisotropic stabilizer and the remaining
  $k$-dense orbits are isolated (cf. case $(RI),(N),(RT)$).


  Now assume that the generic orbit of $X$ is not affine. After
  shrinking $X$, \cref{lemma:maptoG/P} yields an equivariant
  $k$-morphism $\pi:X\to Z$ where $Z$ is a non-trivial complete
  homogeneous $G$-variety. For some $x\in X(k)$ let
  $y:=\pi(x)\in Z(k)$. Then $G_y^{\rm red}$ is a proper parabolic
  subgroup of $G$. Since it has to be conjugate to $P$ (recall
  $\rk_k^{\rm ss}G=1$) we obtain a bijective equivariant morphism
  $G/P\to Z$. Put $\XS:=X\times_ZG/P$ and $\YS:=Y\times_ZG/P$. Then
  $\XS\to X$ is bijective, as well, which, $k$ being perfect, induces
  a bijection $\XS(k)\to X(k)$. Since also the complexities $c(*/P)$
  are unchanged we may replace $X$, $Y$ by $\XS$, $\YS$. Thus we may
  assume that there is an equivariant $k$-morphism $\pi:X\to G/P$.

  Put $y_0:=eP\in G/P$. Then we have $X=G*_PX_0$ where
  $X_0:=\pi^{-1}(y_0)$. Observe that $X_0$ must be irreducible and
  $k$-dense. Moreover, the map $X_0/P\to X/G$ is bijective. Since then
  $c(X_0/P)=c(X/G)\le c_k(X)$ for $Y\subseteq X_0$.

  Let $s\in G(k)$ with $sas^{-1}=a^{-1}$ for all $a\in A$. Then Bruhat
  decomposition of $(G/P)(k)$ (see \eqref{eq:Bruhat}) implies that all
  $k$-points of $X$ are either contained in the closed set $X_0$ or in
  the open set $X_1:=\pi^{-1}(Psy_0)$. Since $Y$ is $k$-dense this
  implies that $Y_1:=Y\cap X_1\ne\leer$. As a $P$-variety, we have
  $X_1=P*_{\Lv}\,sX_0$. Thus $Y_1=\PU\times sY_0$ where
  $Y_0=s^{-1}Y_1\cap X_0$ is $\Lv$-stable and $k$-dense.

  Suppose first that $Y_0\subseteq X_0^A$. Let $x\in Y_0(k)$. Since
  then $A$ is a maximal split torus of $P_x$ all of which are
  conjugate and since $N_P(A)=\Lv$ it follows that $(Px)^A=\Lv x$. By
  continuity, this holds for $x$ in a dense open subset of $Y_0$. From
  this we get
  \[
    c(Y/P)=c(Y_0/\Lv)\le c(X_0/P)=c(X/G)
  \]
  contrary to our assumption. Thus
  $Y_0\cap X_0':=(X_0\setminus X_0^A)\ne\leer$. Now observe that the
  $\Lv$-stabilizer of any $x\in X_0'(k)$ is anisotropic. This implies
  that the set of $x\in Y(k)$ with anisotropic $P$-stabilizer is
  dense. We conclude with \cref{lemma:anisotropic}.

Let  $H$ be the stabilizer of  $x_0\in \pi^{-1}(y_0)$, then $Gx_0\cap \pi^{-1}(y_0)=Px_0\cong P/H$.
Suppose the equality  $c(Y/P)=c_k(X)$ holds. Then also $c_k(X)=c(X/G)$ and  for generic point $x_0\in \pi^{-1}(y_0)$ 
$P$-has a dense orbit in  $Gx_0$, which lies in $Gx_0\cap \pi^{-1}(Psy_0)=PsPx_0\cong PsP/H$.
  By assumptions $H$ has  an open orbit on $PsP/P$. And in particular the adjoint action of $H/H_u$ has an open orbit on $\PU/H_u$ (where $H_u$ is the unipotent radical of $H$).
  If $H$ does not contain a split torus  then  $H/H_u$ is anisotropic and its orbits of the points from ${\PU}/H_u(k)$
  are closed. Combined together this implies the transitivity of the action. Since the point $eH_u$ is fixed, this can happen only for  $H_u=\PU$ and we are in the case of $(U)$-type. 
  
  When $H$ contain a group $A$ and  $H_u\neq \PU$ (otherwise we are in $(U)$-type) the set of $P$-orbits
  contain two orbits with non-anisotropic stabilizers: $Px_0$, and the orbit $Psx_0$ that correspond to the set  $(Px_0)^A=\Lv x$. Since the remaining orbits have anisotropic stabilizer 
  and there is an open orbit with this property, by Corollary \ref{dense_orb_anis} there is only one such orbit. Moreover  when $H_u$ is nontrivial (otherwise we are in type $(RT)$)
   the orbit  $Psx_0$ contains $Px_0$in its closure.
  Indeed these orbits correspond to the H-orbits  $H_u sP/P$ and   $eP/P$ in $G/P$ respectively and $H_u sP/P$ is closed in $PsP/P$. The orbit  $\G_a sP/P$ of  $\G_a\subset H_u$ 
 is  closed in the affine orbit $H_u sP/P$ but not closed in $G/P$.  The closure of $\G_a sP/P$ in $G/P$ gives a $k$-point in the closure of $H_u sP/P$ that is equal to $eP/P$,
 which follows  from Bruhat decomposition (this describes type $(TU)$). 
  
  The calculation of the rank we postpone until next section, where we calculate the lattice. In order to prove $\it (3)$ by Lemma \ref{lemma:anisotropic}
   we have to prove inequality $s_k(Y)>s_k(X)$ in the homogeneously reducible case. In this case for $x\in Y \subset G/H$ the $k$-dense orbit $Px$  maps to $G/P(k)=Psy_0 \cap Py_0$. If it maps to $Psy_0$,
 $P_u$ is acting locally freely and  we apply   \cref{lemma:anisotropic}. Otherwise $Px$ is isomorphic to $P/H$, and  the stabilizer for the action of $\La$ on $P_u\backslash P/H$ is equal to
 $H/\Rad H$ that has  maximal possible dimesion. That implies $s_k(Y)\ge s_k(X)$. 
  \end{proof}

 In the course of the proof of Proposition \ref{prop:main} we proved the following propositions.


\begin{proposition} Let $G$ be semisimple group of semisimple $k$-rank $1$ and
  let $X$ be a $k$-dense homogeneous $G$-variety. Then one of the
  following cases holds:

   \begin{itemize}

  \item[1.] $X$ is $k$-spherical and  homogeneously irreducible over $k$. Then
    $\fB(X)=\{X,Y_1,\ldots,Y_r\}$ with $r\ge0$.

  \item[2.] $X$ is $k$-spherical homogeneously reducible over $k$.

    \begin{itemize}

    \item[2.1] $X$ is horospherical. Then $\fB(X)=\{X,Y\}$

    \item[2.2] $X$ is not horospherical. Then $\fB(X)=\{X,Y_1,Y_2\}$

    \end{itemize}

  \item[3.] $X$ is not $k$-spherical. Then $\fB(X)=\{X\}$.

  \end{itemize}

\end{proposition}

In the following table we summarize the types of raise for the $P$-orbit and we define the action of the group $W_\alpha:=\langle s_\alpha \rangle$.
   We also note that the introduced action is consistent with the action on the character lattices. 

\begin{proposition} \label{H_cases} Let $G$ be semisimple group of $k$-rank $1$ and
   $H$ is the $k$-subgroup.  Then for $H$-orbits on $Y=G/P$ one of the
  following cases holds:

\begin{itemize}

\item[1.]  $H$ spherical, not reductive.

  \begin{itemize}

\item[$(U)$] $H$ is horospherical;\\ $\fB(Y/H)=\{Y>Z\}$; $\rk_kZ=\rk_kY$; $s_k(Z)=s_k(Y)$; $s\circ[Y]=[Z]$.

\item[$(TU)$] $H$ is not horospherical;\\ $\fB(Y/H)=\{Y>Z_1>Z_2\}$; $\rk_k Z_i=\rk_kY-1$; $s_k(Z_1)=s_k(Z_2)$; $s\circ[Y]=[Y]$; $s\circ[Z_1]=[Z_2]$.

  \end{itemize}
  
    \item[2.] $H$ is spherical, reductive

\begin{itemize}
      
\item[$(A)$] $H^\circ$ homogeneously irreducible, anisotropic;\\ $\fB(Y/H)=\{Y\}$; $s\circ[Y]=[Y]$.

\item[$(RT)$] $H$ is homogeneously reducible. $H=L_0.A$ with $L_0\subseteq \La$;\\
  $\fB(Y/H)=\{Y>Z_1,Z_2\}$; $\rk_kZ_i=\rk_kY-1$; $s\circ[Y]=[Y]$, $s\circ[Z_1]=[Z_2]$.

\item[$(RI)$] $H$ is homogeneously  irreducible, isotropic;\\ $\fB(Y/H)=\{Y>Z\}$; $\rk_kZ=\rk_kY-1$; $s\circ[Y]=[Y]$; $s\circ[Z]=[Z]$.
\item[$(N)$] $H$ is homogeneously  irreducible, $H^0$ is homogeneously reducible. $H=L_0.A.\<s\>$ with $s\in N_G(a)\setminus C_G(A)$;\\ $\fB(Y/H)=\{Y>Z\}$; $\rk_kZ=\rk_kY-1$;
$s\circ[Y]=[Y]$; $s\circ[Z]=[Z]$.

\end{itemize}

\item[3.] $H$ is not spherical.

  \end{itemize}

\end{proposition}


Our next aim is to describe how the character lattices of $P$-orbits are related with the restricted Weyl group action in the semi-simple split rank one case.
Let $P_0$ be the maximal normal subgroup of $P$ containing $[P,P]$ such that $\Xi_k(P)=\Xi_k(P/P_0)$. 
If $P_x$ is the stabilizer of $x$ then the calculation of $\Xi_k(P/P_x)$  reduces to study of the $A$-semiinvariant functions
on the  $P_0$-orbit space of $P/P_x$ that reduces to the calculation of the quotient space of $P/P_0\cong A$ by the image of $P_x$ in $A$.

\subsection{Affine irreducible case: types $(RI),(N)$}  

 Let $X_0$ be the open $P$-orbit in $X$ and $y_i\in Y_i$, is the point in the closed orbit with  the stabilizer  equal to $P_H=P\cap H$.
By Lemma \ref{lemma:affine} for $x\in X_0$ any subgroup of $P_x$ which is $k$-isomorphic to $\Bbb G_m$ belong to $Z(G)$. Then the reduced connected component of  the image of $P_x$  in $A$ is $s$-fixed, thus 
   $$s\Xi_k(X)_{\Bbb Q}=\Xi_k(X)_{\Bbb Q}.$$  Let us recall that for  $s$ we can take  an  element of $H$, such that $s\lambda(t)s^{-1}=\lambda(-t)$.  Since $P_x$ contains $\lambda$ and $Z(G)$ is fixed by $s$
   then the image of $P_x$ in $P/P_0$ is fixed by $s$. Then     $$s\Xi_k(Y)_{\Bbb Q}=\Xi_k(Y)_{\Bbb Q}.$$

\subsection{Reducible case: types $(U),(TU),(RT)$ }

After shrinking $X$ we have a $G$-morphism $\pi: X\rightarrow G/P$.  
 
{\bf Type $(U)$.} When $X$ is horospherical i.e. $H \subset P$ and $H_u=\PU$, $X$ consists of the closed $Py:=\pi^{-1}(eP)$ and the open orbit $\pi^{-1}(PsP)$. 
Since $sPsP\supset P$, the orbit $(sPs)y$ contains $Py$ and by the next lemma each orbit $P_uy$ is contained in  $(sP_us)y$. 
We have a map $$k[X]^{(P)}_\chi \rightarrow k[Y]^{(P)}_{s\chi} $$ obtained by applying $s$ to $f$  to get  a $sPs$-semiinvariant function of weight $s\chi$, and  restricting it to $Y$.  
The resulting function $sf|_Y$  is $P_0$-invariant since $sf$ is constant on $sP_0s$-orbits and for each $y\in Y$ the orbit $P_0y$ is contained in  $(sP_0s)y$ that is verified by the following lemma.

\begin{lemma}\label{PUOrb} Let $X$ be a variety with an action of a minimal parabolic group $P$.
Let  be a  $k$-group $H$ group with  anisotropic quotient  $H/H_u$  with the $k$-action  on $X$ such that
 for some $x\in X_k$ we have $Hx\subset Px$. Then we have $H_ux\subset \PU x$ and $Hx\subset P_0x$.
\end{lemma}
\begin{proof}
 Since $\PU$ is a  normal subgroup in $P$  the elementary group $LA$ acts transitively on the set of $\PU$-orbits in $Px$,
 moreover there exists  the geometric quotient  $\pi_{\PU}:Px\rightarrow Z$ by $\PU$ where $Z$ is the homogeneous space for $LA$.
To  prove the lemma we need to show that  $\pi_{\PU}(H_ux)$ is a point.
Since $H_u$ is generated by one dimensional additive $k$-subgroups  $\Bbb G_a$ it is sufficient to prove the lemma for $H_u=\Bbb G_a$, which is the assertion next lemma.

\begin{lemma} \label{A1_G/H}
There are no nontrivial $k$-morphisms from affine line 
$\Bbb A^1_k$ to homogeneous space of an elementary $k$-group $G$.
\end{lemma} 

If $Hx\subset P_0x$ does not hold, we have nontrivial map from $Hx$ to homogeneous space of the split torus $A$ (since $A$ acts transitively of the set of $P_0$-orbits in $Px$).
Since the  $\pi_{\PU}(H_ux)$ is already a point. This map gives a nontrivial map from $H_u\backslash Hx$ which is a homogeneous space of $H/H_u$ 
 to a homogeneous space of the split torus $A$. Which is impossible since by Rosenlicht theorem there are no invertible non-constant $k$-regular functions (which are $k$-characters) on the anisotropic group.
\end{proof}

{\bf  Types $(TU),(RT).$} $X=G/H$ is not horospherical, $H\subset P$ and $H_u\subset P_u$. In this case we can assume that there is a split non-central torus $A_0$ contained in $H$. After conjugation of $A_0$ by $p\in P_k$ 
we can take one-parameter subgroup $\lambda_0\subset A_0$ with $P(\lambda_0)=P$ and such that $M=Z_G(\lambda_0)$
(We note that it may be different from initial $\lambda$ that defines $s$, but these groups have the same centralizer). 
This allows to define the Levi decomposition of $H=\Lv_H H_u$ (where $\Lv_H=Z_H(\lambda_0)$).
   We have $p((G/H)^{A_0})=eP\cup sP$ and  $Y_0=p^{-1}(eP)=P/H\cong \Lv/L_H\times \PU/H_u$ set of  $A_0$-fixed points can be identified with $\Lv/\Lv_H$ since $(\PU/H_u)^{A_0}=e$.
   Also the  $p^{-1}(sP)^{A_0}\cong \Lv/\Lv_H$ is obtained by applying  $s$ to $Y_0^{A_0}$   is contained in $Y_1$ and not in $X_0$. Comparing stabilizers 
   of this points we get $$s\Xi_k(Y_0)_{\Bbb Q}=\Xi_k(Y_1)_{\Bbb Q}.$$
   If $S\subset P$ is one-dimensional non-central torus of $G$ then  $(X_0\cap p^{-1}(sP))^S=\emptyset$ (see Lemma \ref{lemma:anisotropic}). Since  $k$-split tori that are contained in $P_x$ for $x\in X_0$ are central in  $G$,
   then image of the projection of  $P_x$ to $\AXk \cong P/P_0$ is fixed under $s$ and $$s\Xi_k(X)_{\Bbb Q}=\Xi_k(X)_{\Bbb Q}.$$

\section{The general case}\label{s:generalcase}

We continue to assume that $k$ is an infinite perfect field.

\begin{lemma}\label{lemma:finite}

  Let $P$ be a connected linear algebraic $k$-group and let $X$ be a
  $k$-dense $P$-variety. Assume that $c(Y/P)\le c(X/P)$ for all
  $k$-dense $P$-stable closed subvarieties $Y\subseteq X$. Then there
  are at most finitely many $k$-dense $P$-stable closed subvarieties
  $Y$ with $c(Y/P)=c(X/P)$.

\end{lemma}

\begin{proof}

  Suppose there are infinitely many $Y$ with $c(Y/P)=c(X/P)$. Call
  these $Y$ exceptional. Then there is some $d\in\mathbb N$ such that
  there are infinitely many exceptional $Y$ with $\dim Y=d$. Let
  $Z'\subseteq X$ be the closure of their union. Then one of its
  irreducible components, say $Z$, will be the closure of infinitely
  many exceptional \emph{proper} subvarieties. Since $Z$ is a
  $P$-stable $k$-dense subvariety we have by assumption
  $c(Z/P)\le c(X/P)$. Let $Z^0\subseteq Z$ be the open sheet. Then
  every proper $P$-stable subvariety $Y$ of $Z^0$ has
  $c(Y/P)<c(Z/P)\le c(X/P)$. In particular, its closure $\Yq$ is not
  exceptional. Thus $Z$ can't be the closure of a set of exceptional
  subvarieties.
\end{proof}

\begin{theorem}\label{MAIN_COMP}

  Let $G$ be a connected reductive $k$-group with minimal parabolic
  $P\subseteq G$. Let $X$ be a locally linear $k$-dense
  $G$-variety. Then $c(Y/P)\le c_k(X)$ for all $k$-dense $P$-stable
  subvarieties $Y\subseteq X$. Moreover, the set of all closed $Y$
  with $c(Y/P)=c_k(X)$ is finite.

\end{theorem}

\begin{proof}

  Suppose there is $Y$ with $c(Y/P)>c_k(X)$. Then choose one of
  maximal dimension. Clearly $Y\ne X$. Suppose $GY$ is a proper
  subvariety of $X$. Then $c_k(GY)\le c_k(X)$ by [KK] and
  $c(Y/P)\le c_k(GY)$ by induction on $\dim X$. This contradiction
  shows that $GY=X$. Hence, there is a subminimal parabolic $P_\alpha$
  (for some $\alpha\in S_k$) with $Y\ne Z:=P_\alpha Y$. Then
  $c(Z/P)\le c_k(X)$ by maximality.

  Let $U_\alpha$ be the unipotent radical of $P_\alpha$. It is also
  contained and normal in $P$. The group
  $\Gq:=P_\alpha/U_\alpha$ is of semisimple $k$-rank $1$ and contains
  $\Pq:=P/U_\alpha$ as a minimal parabolic. Rosenlicht's theorem
  asserts that $Z$ contains a $P_\alpha$-stable open subset
  $Z^0\subseteq Z$ such that the orbit space $\Zq:=Z^0/U_\alpha$
  exists. Every $P_\alpha$-orbit meets $Y$. So the intersection
  $Y^0:=Y\cap Z^0$ is non-empty and the quotient $\Yq:=Y^0/U_\alpha$
  exists. Clearly, both $\Zq$ and $\Yq$ are $k$-dense. Moreover, $\Gq$
  acts on $\Zq$ and $\Yq$ is a $\Pq$-stable subvariety with
  $c(\Zq/\Pq)=c(Z/P)\le c_k(X)$ and $c(\Yq/\Pq)=c(Y/P)$. Thus, the
  first assertion follows from \cref{prop:main}. The second assertion
  is now implied by \cref{lemma:finite}.
\end{proof}

\begin{corollary}

  Let $G$ be a connected reductive $k$-group with minimal parabolic
  $P\subseteq G$. Let $H\subseteq G$ be a $k$-subgroup and let
  $Z\subseteq Y:=G/P$ be a closed $H$-stable $k$-dense
  subvariety. Then $c(Y/H)\le c(Y/H)=c_k(G/H)$. Equality holds for
  only finitely many $Y$.

\end{corollary}

\begin{corollary}\label{cor:finite}

  A $k$-spherical $G$-variety contains only finitely many $P$-orbits
  $Y$ with $P(k)\ne\leer$.

\end{corollary}

\begin{corollary}

  Let $k$ be a local field of characteristic $0$ and let $X$ a
  $k$-spherical variety. Then $P(k)\backslash X(k)$ is
  finite. Moreover, if $X=G/H$ is homogeneous then
  $P(k)\backslash G(k)/H(k)$ is finite.

\end{corollary}

\begin{proof}

  By \cref{cor:finite}, the set $X(k)$ is covered by finitely many
  $P$-orbits $Y_1,\ldots,Y_n$. Now each intersection
  $Y_i\cap X(k)=Y_i(k)$ decomposes into finitely many $P(k)$-orbit by
  \cite{BorelSerre}*{Cor.\ 6.4}. The second assertion follows from the
  first and $G(k)/H(k)\subseteq X(k)$.
\end{proof}


Finally consider the group $F(W_k)$ defined by  the set of generators $\widetilde{s}_\alpha$ for all  $\alpha \in \Pi_k$ 
(where $ \Pi_k$ is the set of simple roots of the restricted root system $\Delta_k$)
   and by the set of relations $\widetilde{s}_\alpha^2=e$ for all  $\alpha \in \Pi_k$.
Let us introduce the action of  $F(W_k)$ on the set $\fB(X)$ as follows.
For any $\alpha\in \Pi_k$  and $Y\in \fB(X)$ 
the set $P_\alpha Y$ is also a $k$-dense closed subset.
 Consider a space $P_\alpha Y/\Rad_u(P_\alpha)$ that is a rational $P_\alpha$-equivariant  quotient by  $\Rad_u(P_\alpha)$.
  Since the set of $k$-points of the homogeneous space of unipotent group over a perfect field $k$ is always non-empty 
  the preimage induces the inclusion $\fB(P_\alpha Y/\Rad_u(P_\alpha)) \subset \fB(P_\alpha Y)$.
Also by Theorem \ref{MAIN_COMP} $c_k(P_\alpha Y)=c_k(Y)=c_k(X)$, thus $\fB(X) \supset \fB(P_\alpha Y)$.
We have already defined the action of $\widetilde{s}_\alpha$ on $\fB(P_\alpha Y/\Rad_u(P_\alpha))$ according to Proposition \ref{H_cases}, so let us define 
$\widetilde{s}_\alpha\circ [Y]$ in a consistent way.

\begin{conjecture} For a field $k$ such that $\chr k\neq 2$ the action of $F(W_k)$ on the set $\fB(X)$ factors through  the restricted Weyl group $W_k$.
\end{conjecture}

As it was already  mentioned this conjecture is known for algebraically closed fields (cf. \cite{BORB}). 
We prove this conjecture  in the real case in full generality  in Section 5 and  for the action on the  subset $\PS(X)\subset \fB(X)$ of principal families  of $P$-orbits  for ${\rm char} \ k=0$ in Section 6. 

\begin{remark}  Let $X=G/H$ for some $k$-group $H$. Then we can also define the action of $F(W_k)$ on the set $\fB_H(G/P)$ (i.e. on the set of $k$-dense orbits with maximal $H$-complexity) in the following way. 
For each $\alpha\in \Pi_k$ let us fix a morphism $\pi_\alpha:G/P\rightarrow G/P_\alpha$.
For $Hy \subset G/P$ denote by $\widetilde{P}, \widetilde{P}_{\alpha}$ the stabilizers in $G$ of $y$ and $\pi_\alpha(y)$ respectively.
 By the surjectivity of the  map $G(k)\rightarrow (G/P_\alpha)(k)$, the set  $\fB
_H(\pi_{\alpha}^{-1}\pi_{\alpha}(Hy))$ is bijective to the set
$\fB_{H\cap \widetilde{P}_{\alpha}} (\widetilde{P}_{\alpha}/\widetilde{P})$. Consider the image $H\cap \widetilde{P}_{\alpha}$ in $ \widetilde{P}_{\alpha}/\Rad_u( \widetilde{P}_{\alpha})$.
This reduces the situation to a semisimple split rank one case (see Proposition \ref{H_cases}) and
we define  the action of $\widetilde{s_\alpha}$  on $\fB_H(\pi_{\alpha}^{-1}\pi_{\alpha}(Hy))$ to be  consistent  with the action in  the rank one case.
Despite the fact that the  identification of the considered sets $\fB_H$ and $\fB_{H\cap \widetilde{P}_{\alpha}}$ is non-canonical and depends on the choice of $y$ 
a brief look at the cases in Proposition \ref{H_cases} shows that the introduced action does not depend on it.  
As a result we have $F(W_k)$-equivariant embedding of $\fB_H(G/P)$ to $\fB(G/H)$. Moreover $\fB(G/H)=\bigcup_{x\in (G/H)(k)} \fB_{H_x}(G/P)$. 
\end{remark}

\section{$W_k$-action on $\fB(X)$ in the real case}\label{s:realcase}

In \cite{BORB} the spherical case the  $W$-action on the set of  $B$-orbits was introduced by two  different methods via the Hecke algebra and via equivariant
geometry of cotangent bundle that worked only for the orbits of maximal rank. Here we give the way to generalize  the first  method to the $W_k$-action
on the set of $\fB(X)$ in the real  case i.e. $k=\Bbb R$.

For an $H$-variety $X$ consider the category of $(H, {\Gal}(\Bbb C/\Bbb R))$-equivariant  sheaves of $\Bbb F_2$-vector spaces  on $X$ (i.e. on $X(\Bbb C)$) 
that are constructible with respect to the \'etale topology
and denote its  Grothendieck group by $\mathfrak S(X,H)$. Nevertheless, we consider the sheaves as sheaves with respect to the Hausdorff topology. Acoordingly, we use singular cohomology.

We recall that for a $H$-equivariant morphism $f:X\rightarrow Y$ the push forward  $f_!:\mathfrak S(X,H)\to\mathfrak S(Y,H)$ is defined as  $f_![\mathcal F]:=\sum_{i=0}^n (-1)^i[R^if_!\mathcal F]$.
For $\mathcal F_1\in  \mathfrak S (X_1\times X_2,G)$ and $\mathcal F_2\in  \mathfrak S (X_2\times X_3,G)$ we define a product by
$$\mathcal F_1\circ \mathcal F_2=p_{13!}[p^*_{12}\mathcal F_1\otimes p^*_{23}\mathcal F_2]\in \mathfrak S (X_1\times X_3,G).$$

  For $H\subseteq G$ and an $H$-variety $X$ we have an induction functor
  $${\rm ind}^G_H:   \mathfrak S(X,H) \rightarrow \mathfrak S(G*_H X,G): \  \mathcal F \mapsto G*_H\mathcal F=(p_*q^*\mathcal F)^H,$$
for the quotient $p:G\times X\rightarrow G*_HX$ and the action morphism $q:G\times X\rightarrow X$.
 
 When $X$ is already a $G$-variety then there is an isomorphism  $G*_HX\cong G/H\times X$ given by $(g,x)\rightarrow (gH,gx)$.
 Thus, we have an isomorphism $$\mathfrak S(X,H) \cong \mathfrak S(G*_H X,G)\cong \mathfrak S(G/H\times X,G)$$ which provides  $\mathfrak S(X,H)$ with the structure of 
 a $\mathfrak S(G/H,H)\cong \mathfrak S(G/H\times G/H,G)$-module.
 For an inclusion $i:Z\subset X$ by $\Bbb F_{Z}$ we denote the pushforward $i_!\Bbb F_Z$ of a constant sheaf on $Z$.

 The next step is to take into account only those sheaves that have nonzero stalks only on $\Bbb R$-dense $H$-orbits. 
 Let $\mathfrak S(X,H)_{0}\subseteq\mathfrak S(X,H)$ be the subgroup
    generated by all sheaves $\mathcal F$ such that the stalk $\mathcal F_x$ is even dimensional for every $x\in X(\Bbb R)$.
 And let $S(X,H):=\mathfrak S(X,H)/ \mathfrak S(X,H)_{0}$ which is $\Bbb F_2$-module since $2[\mathcal F]=[\mathcal F^{\oplus 2}]$. 

\begin{lemma} Let $\mathcal F\in \mathfrak S(X,H)_0$. Then:
  \begin{itemize}

    \item If $f:X\to Y$ is an $H$-equivariant
    $\Bbb R$-morphism. Then $f_![\mathcal F]\in \mathfrak S(Y,H)_0$.

  \item If $H\subseteq G$ such that If
    $G(\Bbb R)\rightarrow (G/H)(\Bbb R)$ is surjective then
    ${\rm ind}^G_H[\mathcal F]\in \mathfrak S(G*_HX,G)_0$.
  \end{itemize}
\end{lemma}

\begin{proof}
  We mention only the nontrivial steps. First notice that since $G(\Bbb R)\rightarrow (G/H)(\Bbb R)$ is surjective, the condition $Y(\Bbb R)=\emptyset$  for an $H$-invariant $Y\subset X$  is equivalent to
   $(G*_HY)(\Bbb R)=\emptyset$ (Indeed if $y\in(G*_HY)(\Bbb R)$ one can find $g\in G(\Bbb R)$ in such that $gy$ maps to $eH$ under projection to $G/H$ i.e. $gy\in Y(\Bbb R)$).
  This implies that if $\mathcal F\in \mathfrak S(X,H)_{\emptyset}$ then ${\rm ind}^G_H(\mathcal F)\in \mathfrak S(G*_H X,H)_{\emptyset}$.
 Also  an exterior tensor product of a sheaf without $k$-points in its support with any sheaf gives a sheaf without $k$-points in its support.

For any closed  subset $Z\in X$ and the complement $U=X\setminus Z$ there is a long exact sequence in cohomology with compact supports and a
 corresponding long exact sequence for higher direct images:
 $$
 \ldots \rightarrow R_Z^if_{!}\mathcal F  \rightarrow R^if_{!}\mathcal F   \rightarrow  R_U^if_{!}\mathcal F   \rightarrow R^{i+1}_Zf_{!}\mathcal F \rightarrow\ldots
  $$
This allows us to express a push forward in $K$-theory as 
$$
f_![\mathcal F]=\sum_{i=0}^n (-1)^i[R_Z^if_!\mathcal F] +\sum_{i=0}^n (-1)^i[R_U^if_!\mathcal F]. 
$$
Extending this to any stratification $X\supset \ldots \supset Z_i \supset Z_{i+1} \supset \ldots$ by the  closed subsets such that $\mathcal F|_{Z_{i}\setminus Z_{i+1}}$ is  a constant sheaf of rank $r_i$, we get 
$$f_![\mathcal F]=\sum_{j=0}^n (-1)^{j}[R_{Z_i/Z_{i+1}}^jf_!\mathcal F]. $$

 We may choose the stratification in such  a way  that $\mathcal F$ is trivial  on $Z_i \setminus Z_{i+1}$ and  since it has even rank, then
 $R_{Z_i / Z_{i+1}}^if_! \mathcal F$  is also a sheaf of even rank.  Thus   $f_![\mathcal F]\in \mathfrak S(Y,H)_{ev}$.

 \begin{proposition}\label{Smith} Assume that $\Gamma$ is a finite group acting freely on the manifold $X$. And $\mathcal F$ is the $\Gamma$-equivariant locally constant sheaf on $X$. 
 Then $\chi(X,\mathcal F)=|\Gamma|\chi(X/\Gamma,\mathcal F/\Gamma)$. When $|\Gamma|=p$ we also have  
 $\chi(X,\mathcal F)= \chi(X^\Gamma, \mathcal F)\md p$. 
  \end{proposition}

Let  $X$ be a complex manifold  
 such that $X(\Bbb R)= \emptyset$, $\Gamma=\Gal(\Bbb C/\Bbb R)$, which act freely on $X(\Bbb C)$ and   $\mathcal F$-be a constructible sheaf on $X$ defined over $\Bbb R$ (i.e. $\Gal(\Bbb C/\Bbb R)$-equivariant). 
Applying Proposition \ref{Smith} to the restriction of the morphism  $f:X\rightarrow Y$ to the open subset $U\subset Y$ such that $f^{-1}(U)(\Bbb R)=\emptyset$
we get that the rank of $f_![\mathcal F]$ which is the alternating sum of the dimensions of higher direct images is equal to $\chi(f^{-1}(U)(\Bbb C),\mathcal F)=\chi(f^{-1}(U)(\Bbb R),\mathcal F)= 0\md 2$.
\end{proof}

 \begin{corollary}
   If $G(\Bbb R)\rightarrow (G/H)(\Bbb R)$ is surjective, then
   $S(G/H,H)$ is an $\Bbb F_2$-algebra and $S(X,H)$ is an $S(G/H,H)$-module.
    \end{corollary}
 
Let us notice that the $G$-orbits for diagonal action on $G/P\times X$ are in bijective correspondence
 with  $P$-orbits in $X$. Since $G(k)/P(k)=(G/P)(k)$ the $G(k)$ orbits on the set of $k$-points of $G/P\times X$ are in bijective correspondence with $P(k)$-orbits in $X(k)$.
 This identifies  $S (G/P\times X,G)$ with $S (X,P)$. It also provides the group $S(G/P,P)$ with the structure of algebra and the group $S (X,P)$ with the structure of   $S (G/P,P)$-module.
For a group algebra  $\Bbb F_2[W_k]$ we have the following group\footnote{Later we shall prove that 
 it is the isomorphism of algebras.}  isomorphism  $S(G/P,P)\cong \Bbb F_2[W_k]$, so we have a map $\Bbb F_2[F(W_k)]\rightarrow S(G/P,P)$.

 We are going to define the action  group $F(W_k)$ on the set $\fB(X)$  by associating to
 each $Z\subset \fB(X)$ the corresponding element  $\Bbb F_{Z^\circ}\in S(X,P)$ (where $Z^\circ$ is the open $P$-orbit in $Z$) 
 and by assigning to $\widetilde{s}_\alpha$  the action of   $s_\alpha \in S(G/P,P)$ on $S(X,P)$ that is defined 
 as multiplication by the sheaf $\Bbb F_{Ps_\alpha P/P}\in S(G/P,P)$  supported on $Ps_\alpha P/P\subset G/P$.
 This action can be described by the following formula:
 \[s_\alpha\circ [\mathcal F]=\psi_![{\rm ind}_P^{P_\alpha} \mathcal F]- [\mathcal F],   \label{ind}\]
which is a corollary of a Bruhat decomposition of $(P_\alpha/P)(\Bbb R)$ and the following diagram.

$$
\xymatrix{
G/P\times X \ar[d]&\ar[l] P_\alpha/P\times X&\ar[l] \ar[dll]^\psi P_\alpha*_PY&\\
P_\alpha Y& & &\\
}
$$

where the horizontal right map is given by $[g,y]\mapsto (gP,gy)$.

The $P$-equivariant locally constant sheaf with the support on $Pz$ is defined by the irreducible representation of $\rho:P_z/P_z^0\rightarrow GL(V)$, let us denote such sheaf  by $[z,\rho]$.
 Then \ref{ind} can be rewritten in the following way:

$$s_\alpha\circ [z,\rho]= \sum_{i=0}^{2}(-1)^i[P_\alpha\times^H\mathcal H^i]-[z,\rho].$$

We have the following table for the $s_\alpha$-action according to the type of raise.
Let us notice that due to Bruhat decomposition $(P_\alpha/P)(\Bbb R)$ is smooth one point compactification of a disk and thus is diffeomorphic to the sphere $\Bbb S^n$.
In the calculations below we used  the equalities $\chi(\Bbb S^n)=0\md 2$ and
  $\chi(X(\Bbb C),\mathcal F)=\chi(X,\mathcal F)= 0\md 2$.

\begin{itemize}

\item[$(U)$] In this case $Z$ is a point. Thus
$$s_\alpha\circ [y]=(\chi(\Bbb S^n)-\chi(pt))\Bbb F_{\Bbb S^n}  -[y]=[z],$$
$$s_\alpha\circ [z]=\chi(pt)\Bbb F_{\Bbb S^n}  -[z]=[y].$$

\item[$(TU)$]  The closure of $Z_1$ is homeomorphic to  sphere $\Bbb S^k$ and $Z_2=pt$.
$$s_\alpha\circ [y]=(\chi(\Bbb S^n)-\chi(\Bbb S^k))\Bbb F_{\Bbb S^n}-[y]=[y],$$
$$s_\alpha\circ [z_1]=(\chi(\Bbb S^n)-\chi(pt))\Bbb F_{\Bbb S^n}  -[z_1]=[y]+[z_2],$$
$$s_\alpha\circ [z_2]=\chi(pt)\Bbb F_{\Bbb S^n}  -[z_2]=[y]+[z_1].$$
      
\item[$(A)$] Here we have $s_\alpha\circ [y]=[y].$

\item[$(RT)$] Here $Z_i=pt$, thus
  $$s_\alpha\circ [y]=(\chi(\Bbb S^n)-2\chi(pt))\Bbb F_{\Bbb S^n}-[y]=[y],$$
$$s_\alpha\circ [z_1]=\chi(pt)\Bbb F_{\Bbb S^n}  -[z_1]=[y]+[z_2],$$
and the same for the other closed orbit.

\item[$(RI)$] In this case $Z=H/P_H$, and again by Bruhat decomposition $H/P_H\cong \Bbb S^k$.
$$s_\alpha\circ [y]=(\chi(\Bbb S^n)-\chi(\Bbb S^k))\Bbb F_{\Bbb S^n}-[y]=[y],$$
$$s_\alpha\circ [z]=\chi(\Bbb S^k)\Bbb F_{\Bbb S^n}-[z]=[z],$$

\item[$(N)$] Here $Z$ is the union of two points which is actually $\Bbb S^0$.
$$s_\alpha\circ [y]=(\chi(\Bbb S^n)-\chi(\Bbb S^0))\Bbb F_{\Bbb S^n}-[y]=[y],$$
$$s_\alpha\circ [z]=\chi(\Bbb S^0)\Bbb F_{\Bbb S^n}-[z]=[z].$$

\end{itemize}

In order to make out of this action the action of the restricted Weyl group we consider the decreasing filtration  by the dimension of the support of the sheaf
$$S(X,P)_i:=\{\text {generated by } i_*\Bbb F_Z\in S(X,P)|  \text{ for  closed } i:Z\subset X \text{ and }  \dim Z\geq i \}.$$
 An inspection of the above table  shows   the product by generators  $s_\alpha$ of  $S(G/P,P)$
is consistent  with this filtration  (i.e. $s_\alpha\circ S(X,P)_i\subset  S(X,P)_{i-1}$),
thus  ${\rm Gr \, } S(X,P)=\sum_i S(X,P)_i/S(X,P)_{i-1}$ is  the  associated graded $S(G/P,P)$-module 
supplied with the  action of $F(W_k)$ which is consistent with our ad hoc definition from previous section.

Finally   the action of $F(W_k)$  on ${\rm Gr \, } S(X,P)$ 
factors through the action of $W_k$ that follows from the following proposition. 

\begin{proposition} The Hecke algebra $S(G/P,P)$ is isomorphic to a group algebra $\Bbb F_2[W_k]$.
\end{proposition}
\begin{proof}
Let us recall  the $P$-orbits with $k$-points in $G/P$ are parametrized by the elements $w\in W_k$. Also we know that  $P_\alpha\overline{PwP}=\overline{Ps_\alpha wP}$  for $\ell(s_\alpha w)>\ell(w)$
 and  the raise from $P_\alpha\overline{PwP}$ to $\overline{Ps_\alpha wP}$ is of type $(U)$ since all $k$-dense $P$-orbits are of the same $k$-rank. 
 By construction, the action of $s_\alpha\in F(W_k)$ on the element $[PwP]$  is given by $s_\alpha\circ [PwP]=[Ps_\alpha wP]$, which is multiplication by $s_\alpha$. 
 This shows that  $S(G/P,P)$ is identified with $\Bbb F_2[W_k]$.
\end{proof}

Finally we proved. 

\begin{theorem} There is an action of $W_k$ on $\fB(X)$ which factors the action of $F(W_k)$ on $\fB(X)$ and that
is obtained by restriction of the action $S(G/P,P)$ on ${\rm Gr \, } S(X,P)$, where
$W_k$  is realized as a subgroup in $S(G/P,P)$
 by the map on the generators  $s_\alpha \rightarrow [\Bbb F_{Ps_\alpha P/P}]$ and $\fB(X)$ is mapped to the set of lines in  ${\rm Gr \, } S(X,P)$ by  
  assigning  to $Z\in \fB(X)$ the sheaf  $ [\Bbb F_{Z^\circ}]\in  S(X,P)$. 
\end{theorem}

\begin{remark} It is important to notice that  $P_\alpha Y(\Bbb R)$ can decompose as the union of open $P(\Bbb R)$-orbits with different type of the stabilizer. However the $P$-equivariance of the corresponding
constructible sheaves shows that the action of $s_\alpha$  does not depend on the  choice of $P(\Bbb R)$-orbit.
\end{remark}


\section{Action of  $W_k$ on the principal families of $P$-orbits.}\label{s:principal}

From now on let us assume that ${\rm char} \, k=0$. Without loss of generality we can pass to the $G$-invariant subset of smooth points of $X$ and then 
by Sumihiro theorem we can find the covering  $X$ by  $G$-invariant quasi-projective varieties. This allows to assume that $X$ is quasi-projective  so
 by an argument from \cite[Section 5]{inv.mot}  (cf. \cite{KK}) we can pass from $X$ to the total space of some very ample $G$-linearized
line bundle (or just take an affine cone), so we can   assume that $X$ is non-degenerate in the sense of  \cite{inv.mot} or even quasi-affine. 

The aim of this section is  to define the action of the restricted Weyl group on the principal families  of $P$-orbits by means of equivariant geometry of cotangent vector bundle.
Let us recall that for $Y\in \PS(X)$  the action of $\widetilde{s}_\alpha \in F(W_k)$ is defined in the following way:
For  a simple reflection $s_\alpha \in W_k$  consider the corresponding subminimal parabolic subgroup $Q\supset P$ containing preimage of this reflection in $N_G(A)$.
Assume first  that 
 we have $\dim QY>\dim Y$.
If the raise from $Y$ to $Q Y$ is of type $(U)$,  we put  $\widetilde{s}_\alpha\circ Y=Q Y$  otherwise we put $\widetilde{s}_\alpha\circ Y=Y$. 
When $Y=Q Y$,   in the case of existence of $Y'\in \PS
(X)$ such that $Y=QY'$    we put  $\widetilde{s}_\alpha\circ Y=Y'$ (by Theorem \ref{MAIN_COMP} this raise is of type $(U)$),
otherwise we  put $\widetilde{s}_\alpha\circ Y=Y$. 
We are going to prove that  this action of $F(W_k)$  factors through the action of the restricted Weyl group $W_k$.

Let us denote by $\muG$ the composition map $\CT{X}\rightarrow
\gm^*\rightarrow \gm^*/\!\!/G\cong \tm^*/N_G(T)$,
where the last equality is the Chevalley isomorphism,
and consider the morphism:
$$
\CTt{X}:=\CT{X}\times_{\tm^*/W} \tm^*\rightarrow \gm^* \times_{\tm^*/W} \tm^*,
$$
that coincides with the moment map on the first factors and with
identity morphism on the second factors.
For algebraically closed field we have the following theorem of F.Knop.

\begin{lemma}{\cite{BORB}}\label{MAX_RC} All irreducible components of $\CTt{X}:=\CT{X}\times_{\tm^*/W} \tm^*$
have the same dimension. They map onto $\CT{X}$ and $W$ acts
transitively on them via its natural action on
the second factor.
\end{lemma}

For algebraically non-closed fields these results can be refined in the following way.
Consider the following sequence of finite maps:
$$\tm^*\stackrel{/W_L}{\longrightarrow}  \lv^*/\!\!/L \stackrel{/W_k}{\longrightarrow}\lv^*/\!\!/ N_G(L) \longrightarrow \gm^*/\!\!/G\cong \tm^*/W,$$
and let us take  into account that $N_G(L)=N_G(A)$,  $ \lv^*=(\g^*)^A$ and that the map $\lv^*/\!\!/ N_G(L) \longrightarrow \gm^*/\!\!/G$ 
is finite by the result of Luna (see \cite[Thm.6.16]{inv}).

There are a least two candidates for the substitute of  $\CTt{X}$ in an algebraically non-closed case:
$$
\CTt{X,P}:=\CT{X}\times_{\tm^*/W}\lv^*/\!\!/L, \ \ \ \ \ \   \CTt{X,P'}:=\CT{X}\times_{\tm^*/W}\ab^*.
$$
The first one is better related with geometry of cotangent bundle. And the second is important being related with the sections of conormal bundle to $\PU$-orbits which are defined 
by the semi-invariant functions.
These objects are related by the following sequence of maps, where the right down vertical arrow is induced by $L$-equivariant splitting $\ab^*\subset \lv^*$. 
$$
\xymatrix{
& \CTt{X,P'}=\CT{X}\times_{\tm^*/W}\ab^*   \ar[dl] \ar[d]& &    \\
\CT{X}\times_{\tm^*/W}\tm^*      \ar[r]^{/W_L}  &  \CT{X}\times_{\tm^*/W}\lv^*/\!\!/L  \ar@<1ex>[u]  \ar[r]^{/W_k} & \CT{X}\times_{\tm^*/W}\lv^*/\!\!/ N_G(L)  \ar[r]    &\CT{X} 
}
$$
From 
 \cite{BORB} we know that families of $B$-orbits of maximal rank and complexity defined over $K$ (denoted by $\B_{00}(X)$) correspond to connected components of $\CTt{X}$.
Our aim is to associate to each element of $\fB_0
(X)$ 
 the corresponding irreducible  components of $\CTt{X,P}$ 
   and $\CTt{X,P'}$ 
and to prove that the  action of $W_k$  on the irreducible components of $\CTt{X,P}$ 
and $\CTt{X,P'}$ via its action on the right factor
correspond to the action of $F(W_k)$ on  $\fB_0
(X)$ introduced above.
 We have to consider both varieties  $\CTt{X,P}$ and $\CTt{X,P'}$ simultaneously  since it is not a priory  clear that the natural maps between them are 
  inclusions on the irreducible components which is due to the fact that  $\CTt{X,P'}\rightarrow \CT{X}$ is not surjective. The first variety is important because it is more closely  related with  the geometry of the moment map. 
 The second one is interesting since it is related with the geometry of the sections for $P_u$-orbits.


To study the components of $\CTt{X,P}$ let us consider the composition of maps: $$\CT{X}\times_{\tm^*/W}\lv^*/\!\!/L\rightarrow \gm^* \times_{\tm^*/W} \lv^*/\!\!/L\rightarrow \gm^*$$
Let us notice that the last map 
admits a section over $\p_u^{\perp}$ defined
by
$$\tau: \p_u^{\perp} \rightarrow \gm^* \times_{\tm^*/W} \lv^*:    \ \
\lambda \mapsto (\lambda,\lambda |_{\lv^*}).$$
This formula is well defined by the following lemma.

\begin{lemma}\label{adapt} Let $\p \rightarrow \lv$ be an $L$-equivariant projection to the Levi subgroup. Then for  $\xi \in \lv$ the fiber $\xi+\p_u$ is contained in the fiber 
of quotient map $\g^*\rightarrow \g^*/\!\!/G$. Moreover if $\xi_0\in\xi+\p_u$ is semisimple then
$\xi_0\in P_u\xi$.  
\end{lemma}
\begin{proof} Let $\lambda(t)\in \Ch(A)$  be a one-parameter subgroup
of $A$ adapted to $P$. Then we have $\lim_{t\rightarrow 0} \lambda(t)\xi_0=\xi$ for all $\xi_0\in\xi+\p_u$, that implies the first claim.
If $\xi_0\in\xi+\p_u$ is semisimple,  then there exists an element $p_u\in P_u$ such that $\Ad(p_u)\xi_0\in \lv.$ Since $\lv$ is fixed by $\lambda$ and 
$\lim_{t\rightarrow 0}\lambda(t)p_u\lambda(t)^{-1}=e$ we have
$
\Ad(p_u)\xi_0=\lim_{t\rightarrow 0}\Ad (\lambda(t)p_u)\xi_0=\Ad(\lim_{t\rightarrow 0}\lambda(t)p_u\lambda(t)^{-1})\lim_{t\rightarrow 0}\Ad(\lambda(t))\xi_0=\xi.
$\end{proof}
This gives us the embedding of $\Mom^{-1}(\p_u^{\perp})$ in
$\CTt{X,P}$  
via $\widetilde{\tau}:\eta  \mapsto (\eta,\Mom(\eta)|_{\lv^*}).$
The set $X$ also defines the distinguished component of $\CT{X}\times_{\tm^*/W}\lv^*/\!\!/ N_G(L)$ which is denoted by $\CT{X}$. 
The main result of this section is the following.

\begin{theorem}\label{Waction}  There exists the  action of $W_k$ on the set $\fB_0
(X)$ of the   principal families of $P$-orbits,
and $\fB_0(X)$  is isomorphic to the set  of those $k$-irreducible components of $\CTt{X,P}$ (resp. $\CTt{X,P'}$)
that map to the distinguished component of $\CT{X}\times_{\tm^*/W}\lv^*/\!\!/ N_G(L)$.
 The $W_k$-action on this set of components is  induced by the action on the second factor of $\CT{X}\times_{\tm^*/W}\lv^*/\!\!/L$.
\end{theorem}

The variety $\Mom^{-1}(\p_u^{\perp})$ is the union of the conormal
bundles to all $P_u$-orbits in $X$. 
For an $H$-invariant subvariety $Y$ by $\CN[Y]{H}$ (or briefly by $\mathcal N(Y)$ when $H=P_u$) we denote the conormal bundle to the family of generic $H$-orbits in $Y$.
By $\NT(Y,P)$ (resp. $\NT(Y,P')$)  we denote the image of $\mathcal N(Y,P_u)$ (resp. $\mathcal N(Y,P')$) in  $\CTt{X,P}$ (resp. $\CTt{X,P'}$).
Since each element of $\g^*$ is $G$-conjugate to the element in $\p$, we get that for each component of $\CTt{X,P}$ there is an irreducible component $\widetilde{\mathcal{N}}$ of $\Mom^{-1}(\p_u^{\perp})$ whose $G$-span is dense. Denoting by $Y$ the image of projection of $\widetilde{\mathcal{N}}$ to $X$ we get $\widetilde{\mathcal{N}}=\mathcal{N}(Y)$.

Note  that  the dimension  of  the family of  $P_u$-orbits is equal to $c(Y/P)+\rk_k(Y)+s_k(Y)$ for every $P$-invariant $k$-dense subvariety $Y$. 
The dimension\footnote{
For the split case the dimension of the component of conormal bundle is equal to $\dim X+c(Y/P)+\rk_k(Y)$.} for the corresponding conormal bundle to the foliation is equal to $\dim X+c(Y/P)+\rk_k(Y)+s_k(Y)$,
which is the same for all  principal families.

\begin{definition}
The elementary radical $\Rad_\el H$ (anisotropic radical  $\Rad_\an H$) is the smallest normal subgroup such that  $H/\Rad_\el H$ ($H/\Rad_\an H$) is elementary (anisotropic) subgroup.
\end{definition}

\begin{definition}
  Let $H$ be a connected $k$-group acting on a $k$-dense variety
  $X$. Then the action 
  is called \emph{elementary} 
  (\emph{anisotropic}) if the elementary radical $\Rad_\el H$ (the
  anisotropic radical $\Rad_\an H$)  acts trivially on $X$.
\end{definition}

\begin{remark}
It is easy to prove that $\Rad_\el H$ is generated by $k$-defined unipotent subgroups of $H$, and  $\Rad_\an H$ is generated 
by $\Rad_\el H$ and  $k$-defined subgroups $\Bbb G_m$ of $H$.
\end{remark}

Let us recall the definition of the normalizer of generic $P$-orbit in $X$, which is a parabolic subgroup of $G$.
$$\PXk:=\{ g\in G \ | \ gPx=Px \text{\ for \  } x \text{\ in a dense open subset of } X \}$$

 When we work over the fixed field  $k$ we use the simplified notation $\PX$.
We shall need the following version of the Local Structure Theorem proved by F.Knop and B.Krotz \cite[Prop. 4.6]{KK}.

\begin{theorem}[Generic Structure Theorem]\label{LSTgeneric}

  Let $X$ be a $k$-dense $G$-variety and let $\PX=\LX \PUX$ be the Levi
  decomposition in which $\LX$ is normalized by $A$. Then there exists a smooth affine $\LX$-subvariety
  $X_\el\subseteq X$ such that

  \begin{enumerate}

  \item the action of $\LX$ on $X_\el$ is elementary, all orbits are
    closed, and the categorical quotient $X_\el\to X_\el /\!\!/ \LX$ is a
    locally trivial fiber bundle in the etale topology. In particular there exists the stabilizer $\LXS$ of general position for this action.

  \item the natural morphism $\PUX\times X_\el=\PX*_{\LX} X_\el\to X$ is
    an open embedding.

  \end{enumerate}

  The slice $X_\el$ is unique up to a unique $\LX$-equivariant birational
  isomorphism. More precisely, its field of rational functions can be
  computed as
  \[\label{eq:kXN}
    k(X_\el)=k(X)^{\PUX}=k(X)^{P_u}.
  \]

 \end{theorem}

By $\LXo,\TXo,\AXo$ let us denote subgroups of $\LX,T,A$ which are the kernels  of  the characters of all $P$-semi-invariant functions on $X$. 
By $\ax$ let us denote the orthocomplement to $\axo$ in $\ab$. 

\begin{remark} 
The above theorem implies the existence for the stabilizer of general position for the action $P$ on $X$ as well as existence of stabilizer of general position 
  $L(X):=\LXS\cap L$  for the action of $L$ on   $X_\an$. 
 We also have   $s_k(X)=\dim L/L(X)A$.  
\end{remark}



\begin{lemma}\label{stab} Let $Y\in \PS(X)$. 
Then there exists a sequence $P_{1},\ldots, P_{l}$ of subminimal parabolic subgroups
containing $P$  such that in the sequence $Y_i:=P_i\ldots P_1Y$   the raise from $Y_{i-1}$  to $Y_i$ is of type $(U)$ and $P_l\ldots P_1Y=X$. Moreover there exists a stabilizer of general position
for the action $P$ on $Y_i$ and there exists a stabilizer of general position $L(Y)$ for the action of $L$  on $Y_i\ds P_u$ which is conjugate to $L(X)$ by the action of $w:=s_l\ldots s_{i-1}\in N_G(A)_k$.  
\end{lemma}
\begin{proof}

Let us argue by  decreasing induction on the dimension of $Y_i$. The base of induction is given by the local structure theorem \ref{LSTgeneric}.
Given $Y\in \fB_0(X)$,  consider a  parabolic minimal subgroup  $P_1$ containing $P$ for  which $P_1Y\neq Y$ (Such $P_1$ exists since $G$ is generated by supminimal parabolic subgroups and 
according to \cite{KK} if $\overline{GY}\neq X$ we have $c_k(Y)<c_k(X)$). 
Since the $(c_k(Y),\rk_k(Y))=(c_k(X),\rk_k(X))$ the raise from $Y$ to $P_1Y$ is of type $U$. Moreover there is a $P_1$-equivariant map $\varphi:P_1y\rightarrow P_1/P$.
By induction on the dimension we can assume that the assertion of the lemma is true for $P_1Y$ and $P_l\ldots P_1Y=X$. 
Let us recall that for a general $y\in Y$  the set of $k$-dense  orbits of  $P_1 y$ consists of the open orbit $Ps_1 Py=Ps_1 y=\varphi^{-1}(Ps_1Py_0)$ and the  orbit $Py=\varphi^{-1}(y_0)$ in its closure.
Meantime applying $s_1$ to   $P_1 y$ we get that the same decomposition in the union of   $s_1Ps_1$-orbits. The orbit  $s_1Py$ is closed and  $s_1Ps_1Py$ is open so it is also equal to
$s_1Ps_1y$, which follows from $y_0\in s_1Ps_1Py_0$ and $\varphi^{-1}(s_1Ps_1Py_0)= \varphi^{-1}(s_1Ps_1y_0)$. Let  us also notice that  
$s_1Ps_1Py$ contains $Py$.



Since $(s_1Ps_1)y$ contains $Py$ by  lemma \ref{PUOrb}  the orbit $P_uy$ is contained in  $(s_1P_us_1)y$. In particular if $Z$ is the $L$-invariant section for $P_u$-orbits in $Y$  (that always exists),
then it is also the section for the family of $s_1P_us_1$ in $P_1Y$ and then $s_1Z$ is the section for the family of $P_u$-orbits in $P_1Y$.
Taking into account that $s_1\in N_G(S)_k$ normalizes $L$, this proves that stabilizer of general position for the action of $L$ on $Y_i\ds P_u$ is equal to $L(Y)=s_1L(P_1Y)s_1=wL(X)w^{-1}$. 
\end{proof}

\begin{corollary}{\label{restr_conorm}} Let $P_1$ be a subminimal $k$-parabolic subgroup properly containing $P$, $s\in N_{P_1}(A)$ be a representative of a simple
reflection. Let $Y\in \fB_0(X)$ be such that $P_1 Y$ be the span of type $(U)$. Then 
 for all $y\in Y$ we have  $$s\CN[P_u]{P_1 y}=\CN[\Ual]{\Bal y}\subset (sPs\cap P_1)\CN[P_u]{Py},$$

 \end{corollary}
\begin{proof}
Since $Py\subset \Bal y$ for a general point of $y\in Y$, the Lemma \ref{PUOrb}  implies that  $P_uy\subset \Ual y$.
Thus the restriction of $\CN[\Ual]{\Bal y} $ to $Py$ is contained in
 $\CN[P_u]{Py}$ and we have  $s\CN[P_u]{P_1 y}=\CN[\Ual]{\Bal y}\subset (sPs\cap P_1)\CN[P_u]{Py}$. 
\end{proof}

\begin{proposition}\label{dim_comp} There is a bijection $Y \mapsto \CN{Y,P_u}$ between
$Y\in \fB_0(X)$ and the $k$-dense irreducible components of
$\Mom^{-1}(\p_u^{\perp})$  that dominate $\mu(\CT{X})$  and have dimension 
$\dim X+c(X/P)+\rk_k(X)+s_k(X)$ which is maximal for the components with this property.
\end{proposition}

The following proposition whose proof is postponed until the end of the section finishes the proof of Theorem \ref{Waction}.

\begin{proposition}\label{Wprop} The action of the group $F(W_k)$ on the set of $\fB_0(X)$ factors through $W_k$
 and is isomorphic to the action on the set  of the $k$-irreducible components of $\CTt{X,P}$ (resp. $\CTt{X,P'}$)
that map to the distinguished component of $\CT{X}\times_{\tm^*/W}\lv^*/\!\!/ N_G(L)$. Where the bijection on the sets is given
by  $Y \mapsto G\CNt{Y,P_u}.$
\end{proposition}

By \cite{KK} for an quasi-affine $X$ there exists a $\chi\in\Ch_k(X)$  that satisfies $\langle \chi,\alpha \rangle\neq 0$ for any $\alpha \in \Delta_{\pux}$.
Such $\chi$ can also be characterized by the property $s_\alpha \chi \neq \chi$ for $s_\alpha \in N_G(A) \setminus \LX$.
We consider   $\Ch_k(X)$ as a lattice in $\ax^*$. The collection of such $\chi \in \ax^*$ which satisfies the above condition is Zariski dense in $\ax^*$ and is denoted by $\ax^{pr}$. 
By $A_0$ let us denote the common kernel in $A$ of characters of all $P$-semi-invariant functions ($\ab=\ab_0+\ax$).
Consider a $\PX$-semi-invariant function $f_\chi\in k(X)^{(Q)}_\chi$,
 following Knop (\cite{inv.mot}) we can define a $\PX$-equivariant map
$$\logd[\chi] \colon X \setminus \dvr (f_\chi) \longrightarrow \mathfrak g^*, \ x\mapsto l_x,
 \ \  where \ \ l_x(\xi)=\frac{\xi f_\chi}{f_\chi}(x).$$

The map
$ \Logd[\chi]:\X\rightarrow \CT{X}$, where $\X=X\setminus \dvr(f_\chi)$, that maps $y\in \X$ to the value of the section 
$f_\chi^{-1}(df_\chi)$ in the point $y$, produces the following commutative diagram:

$$
\xymatrix{
 &
  \CT{X}  \ar[dr]^{\Mom}& & \\
\X \ar[rr]^{\logd[\chi]}  \ar[ur]^{\Logd[\chi]} & &\gm^* &
}
$$

Let us notice that  $\Logd[\chi](X)$ is the subset of $\CN{X}$, since the differential form $f_\chi^{-1}(df_\chi)$ annihilated on the action vector field of $\p_u$.
We also recall  from \cite[(4.5)]{KK}  that the image  $\logd[\chi](\X)$ is equal to $\xi+\pux$.
The variety  $\K_\chi=\psi_\chi^{-1}(\chi)\subset \X$ (or simply $\K$) provides the section for $P_u$-orbits
 (which are also $\PUX$-orbits)  in $\X$. The $\PX$-equivariant map $\Logd[\chi]$ embeds $\K$   in $\CN{X}$,  we put  $\widetilde{\K}:=\Logd[\chi](\K)$.

Before stating next proposition recall that by reductivity of $L$  we have an  isomorphism $\lv\cong \lv^*$ and a splitting $\g^*\cong \lv^*\oplus \lv^\bot$ which are both $L$-equivariant.

\begin{proposition} \label{wZ_section} Let $G$ be $k$-group acting on $k$-dense quasiaffine variety $X$. For $X_\el$ which is birational to $X\ds P_u$,
 consider  $\Mom[\lv]: \CT{X,\el} \rightarrow \lv^*$ the corresponding moment map. For $Y\in \fB_0(X)$
and for the presentation of the element  $w:=s_1\ldots s_l\in W$   consider
$Y_i:=P_i\ldots P_1Y$ and assume that each the raise from    $Y_i$ to $Y_{i+1}$ is of type $(U)$. 
Then

(i) $G\CN{Y}$ is dense in  $\CT{X}$. 

(ii) $\overline{\Mom(\CN{Y,P_u})}=Pw\Mom[\lv](\CT{X,\el})$ and $\overline{\Mom( \CN[P']{Y})}=Pw\ax$.

(iii) $w\K$ is the $L$-invariant rational section for generic $P_u$-orbits in $Y$.



\end{proposition}


\begin{proof} 
For the action of $\LX$ on $X_\el$ by Luna slice theorem there exists a slice $Z$ (which we can assume is defined over $k$), so we have an excellent map:
$\LX *_{\LXS}Z\rightarrow X_\el$, and the corresponding map $\PX *_{\LXS}Z\rightarrow \X$. 
The  $\LX$-orbits on $X_\el$ coincide with the $L$-orbits.   So we have 
$$L*_{\LS }Z\rightarrow X_\el, \ \ \ \ \  \PUX \times (L*_{\LS }Z)\rightarrow \X,$$ 
which gives the following $L$-equivariant etale map to $\CT{X}$:
$$\PUX \times \pux \times (L*_{\LS }\ls  ^\bot) \times \CT{Z}.$$

To embed $\CT{Z}$ into $\CT{Z}$ we used  an isomorphism $\CN{X,P_X}|_Z\iso \CT{Z}$, obtained by  the restriction of the elements of $\CN{X,P_X}\subset \CT{Z}$ to $T_Z$.
Also we have  $$ \CT{L*_{\LS }Z}\isol \CN{X,\PUX}|_{L*_{\LS }Z}\subset \CT{X}.$$

The action of $\Rad_\el \LX$ 
 is trivial on $X_\el$. After taking quotient by $\PUX$ we obtain that $P/\PUX$ is minimal parabolic for $\LX$. 
Since anisotropic group $\LX/\Rad_\an \LX$ does not contain  unipotent elements, the image of    $P$ in   $\LX/\Rad_\an \LX$ coincide with the image of $\La$
and this restriction map is surjective.
  We denote by $\lvo  $ the lie algebra  which is the kernel of this map. 
The lie algebra which is the direct summand  of $\lv$ and which coincide with  the image of $\lv$ in the lie algebra of $\LX/\Rad_\an \LX$ is denoted by $\lve$.

The $P_u$-orbits in $X$  are also normalized by $P_u(L\cap \Rad_\el \LX)A_0$.
Then the conormal bundle $\CN{X}$ to $P_u$-orbits in $X$ is also conormal bundle to   $P_u(L\cap \Rad_\el \LX)A_0$-orbits and its image under the moment map belong 
 to $(\p_u+\lvo  +\ab_0)^\bot=\ax+\lve+\pux$.  

Let us notice that the restriction of the moment map to the subset $\CN{X}$ and subalgebra $\lv$ can be obtained by composition of the moment map    $\CN{X}\rightarrow \lve+\pux$
and the projection to $\lve$. Recall that we have the etale map 
$$\PUX \times (L*_{\LS}\ls  ^\bot) \times \CT{Z}\rightarrow \CN{X,P_u}.$$
In particular $\dim \CN{X,P_u}=\dim \CT{X}-\dim \PUX$
Let $\tm_0^\bot$  be the orthocomplement in $\tm$ of the Lie algebra $\tm_0$ that annihilates $B$-semi-invariant rational $K$-functions on $X$. Then by \cite{inv.mot} the closure of the
image of the moment map $\Mom[\lv]: \CT{X,\el} \rightarrow \lv$ is equal to $L\tm_0^\bot$. 
Taking into account that   $\langle \chi,\alpha \rangle\neq 0$ for $\chi \in \tm_0^\bot$ and $\alpha \in \Delta_{\pux}$ we get  $\PUX\chi=\chi+\pux$.
 So the closure of the image of the map $\CN{X,P_u}\rightarrow\g^*$ contains $\PX\tm_0^\bot=L(\tm_0^\bot+\pux)=L\tm_0^\bot+\pux$.  It cannot be larger
 since  $\Mom[\lv](\CN{X,P_u})=L\tm_0^\bot$ which is the   projection of  $\Mom(\CN{X,P_u})$  to $\lv$.

Since $\PUX^- \chi=\chi+\pux^-$ is transversal to $L\tm_0^\bot+\pux$ in the point $\chi$, then it is also true for generic point $L\tm_0^\bot+\pux$. By looking at differential of $d\Mom$ we get 
that $\PUX^- y$ is transversal to $\CN{X,P_u}$  in the  generic point $y\in \CN{X,P_u}$. This implies that $\dim \PUX^- \CN{X,P_u}=\dim \PUX^-+ \dim \CN{X,P_u}=\dim \CT{X}$ and that $\PUX^- \CN{X,P_u}$ is dense in $\CT{X}$.

 
 To prove $(i)$ and $(iii)$ let us argue by descending induction on the dimension of $Y$.  Assume that $G\CN{P_1 Y,P_u}$ is dense in $\CT{X}$. By corollary \ref{restr_conorm} we have
 $s\CN{P_1 Y,P_u}\subset (sPs\cap P_1)\CN{Y,P_u},$  that implies the density of $G\mathcal N(Y,P_u)$ in $\CT{X}$.

The construction of $\K$ gives the base of the induction for the proof of $(iii)$. 
By taking a rational quotient and passing to the open subset of $Y$ we can assume that existence of a 
 quotient morphism from $Y_1:=P_1  Y$ to the smooth variety $Y_1/(P_1)_u$.
For the conormal bundle  $\CN[(P_1)_u]{Y}$ to the foliation of the $(P_1)_u$-orbits in $Y_1$ 
consider  a natural map to $\CT{Y_1/(P_1)_u}$, which is a composition of natural projection to $\CT{Y_1}$ 
and the quotient by $(P_1)_u$. 
We note that the conormal bundles    $\CN[P_u]{Y_1}$ and $\CN[P_u]{Y}$ to the foliations of  $P_u$-orbits are preimages under this map of the conormal bundles  to the foliations of  the $\overline{P}_u$-orbits in 
$Y_1/(P_1)_u$ and $Y/(P_1)_u$. 
We have a following commutative diagram:

$$
\xymatrix{
\CN[(P_1)_u]{Y_1} \ar[d] \ar[r] & (\p_{1})_u^\bot \subset \g^*  \ar[d] 
  & \\
\CT{Y_1/(P_1)_u}\ar[r]  &   \m^* &
}
$$

 By induction assumption we have the section $\K_1=s_1w\K$ for the set  of general $P_u$-orbits in $Y_1$ and its lifting
 $\widetilde{\K}_1 \subset \CN[P_u]{Y_1}$  such that $ \Mom(\widetilde{\K}_1)=\chi\in \Chpr{Y_1}$.
We have to prove that $s_1 \K_1$
is the section for the $P_u$-orbits in $Y$. From the proof of the Corollary \ref{restr_conorm} we know that $s_1 \K_1$ is the section for 
$s_1P_us_1$-orbits in $Y_1$ and a each generic orbit contain a unique corresponding $P_u$-orbit of $Y$.
 Suppose that $s_1 \K_1$ is not a section for  the set of $P_u$-orbits in $Y$, so we have $z\in \K_1$ 
 (and its lifting $\widetilde{z}\in \widetilde{\K}_1$) such that $s_1 z \notin Y$.
Let us take $u^{-}\in (P_u)^-\cap P_1$ such that $u^{-} s_\alpha z \in Y$.
 The set of $P_u$-orbits in $Y_1$  is not normalized  by $(P_u)^-\cap P_1$.

\begin{lemma}\label{degen_weight} The general weight $\chi_1=s_1w\chi$ of  $\Ch(Y_1)$ is non-degenerate with respect to $P_1$.
\end{lemma}
\begin{proof}
 The set of $P$-orbits in $Y_1$  is not stabilized by $P_1$ and 
in particular not stabilized by $s_1$. Then the points of the section $\K_1=s_1w\K$ are not stabilized by $s_1$. Recalling that the stabilizer of $\K_1$ is equal to $s_1w\LX w^{-1}s_1$ 
which is the centralizer of $s_1w\chi$ thus $s_1$ do not fix it.
\end{proof}


Since   $s_1\chi_1 \neq \chi_1$,
   we have $\Mom(u^{-} s_1 \widetilde{z})=s_1 \chi_1 + \eta_{-}$, where $\eta\in (\p_u)^-\cap \p_1$ and $\eta\neq 0$.
  We have the inclusion $s_1 \widetilde{\K}_\alpha\subset \CN[\Ual]{Y_1}$.
 Moreover since $u^{-} s_1 z \in Y$ the element $u^{-} s_1 \widetilde{z}$ belongs to the restriction of $\CN[\Ual]{Y_1}$ to $Y$.
By the Corollary \ref{restr_conorm} we have $u^{-} s_1  \widetilde{z}\in \CN[U]{Y}$ and in particular $\Mom (u_{-\alpha} s_\alpha  \widetilde{z})\in \p$.
We get a contradiction since $s_\alpha \chi + \eta_{-}\notin \p$, thus we get $u^{-}=e$ and $s_1 \K_1\subset Y$.
As a corollary we also get that $\Mom(\CN{Y})\supset P\aY{Y}$. 
Moreover the points of $w\K$ are stabilized by $w(\LS A_0\Rad_\an \PX) w^{-1}$,

\begin{lemma}\label{not_dense} Let $Y\in \fB(X)$ be such that the raise from  $Y$ to $Y_1:=P_1Y$ is not of type $(U)$. 
Then  the action of $s_1\in W_{P_1}$ on the component of $\CTt{X,P}$ (resp. $\CTt{X,P'}$) corresponding to $\CN[P_u]{Y_1}$ (resp. $\CN[P']{Y_1}$) is trivial.
\end{lemma}
\begin{proof}
 By the overline let us denote the corresponding images of the groups in $L_1\cong P_1/(P_1)_u$.
 In the following diagram  we have to prove 
 that for $L_1$-variety $\overline{Y}_1:= Y_1/(P_1)_u$ the variety $\CT{\overline{Y}_1}\times_{\tm^*/W_1} \tm^*$ (resp. $\CT{\overline{Y}_1}\times_{\tm^*/W_1} \s^*$ ) is irreducible. 
$$
\xymatrix{
	\CT{\overline{Y}_1}\times_{\tm^*/W_1} \s^*    \ar[d] & \ar[l] P_1\CN[P']{Y_1}  \ar[r]     \ar[d] &  \ar[d] \CT{X}\times_{\tm^*/W} \s^* \\
	\CT{\overline{Y}_1}\times_{\tm^*/W_1} \lv^*/\!\!/L   & \ar[l] \CN[(P_1)_u]{Y_1}  \ar[r]    &\CT{X}\times_{\tm^*/W} \lv^*/\!\!/L
}
$$
By the assumption on the raise the triple $(c_{\overline{P}},{\rk}_{\overline{P}},s_{\overline{P}})$  takes strictly lower value on  $\overline{P}$-invariant subvariety of $\overline{Y}_1$ which we denote by $\overline{Y}$.
In case  ${\rk}_{\overline{P}}(Y)<{\rk}_{\overline{P}}(\overline{Y}_1)$ we get that $\Mom (\CN[\overline{P}_u]{{Y}})\subset \ab_Y+\p'$ so the closure of  $G\CN[\overline{P}_u]{\overline{Y}}$
is not equal to $\CT{\overline{Y}_1}$ (by the similar argument  $G\tilde{\tau}\CN[\overline{P}']{\overline{Y}}\subset \CT{X}\times_{\tm^*/W} \ab_Y$ does not provide a component of $\CT{X}\times_{\tm^*/W} \s^*$).

Otherwise the sum $c_{\overline{P}}+{\rk}_{\overline{P}}+s_{\overline{P}}$ (resp. $c_{\overline{P}}+\rk_{\overline{P}}$) takes strictly lower value on  
 $\overline{Y}$   that gives  $\dim\CN[\overline{P}_u]{\overline{Y}}<\dim\CN[\overline{P}_u]{\overline{Y}_1}$ (resp.  $\dim\CN[\overline{P}']{\overline{Y}}<\dim\CN[\overline{P}']{\overline{Y}_1}$). 
Taking into account that $\dim\CN[\overline{P}_u]{\overline{Y}_1}=2\dim \overline{Y}_1-\dim \overline{P}_u$, we get
$$\dim(L_1\CN[\overline{P}_u]{\overline{Y}})=\dim(\overline{P}_u^-\CN[\overline{P}_u]{\overline{Y}})<\dim T^*\overline{Y}_1.$$
(resp. $\dim(L_1\CN[\overline{P}']{\overline{Y}})=\dim(\overline{P}_u^-\CN[\overline{P}']{\overline{Y}})<\dim(\overline{P}_u^-\CN[\overline{P}']{\overline{Y_1}})$
where the last inequality is strict. 
Indeed $\dim(\overline{P}_u^-\CN[\overline{P}']{\overline{Y}})=\dim \overline{P}_u+\dim \CN[\overline{P}']{\overline{Y}})$ since $\overline{P}_u^-$-orbits of the generic points in  $\CN[\overline{P}']{\overline{Y_1}})$ are transversal to this variety, that was  checked in the beginning of the proof of Proposition \ref{wZ_section}).

Thus $L_1\CN[\overline{P}_u]{\overline{Y}}$ is not dense in $T^*\overline{Y}_1$ (resp. $L_1\CN[\overline{P}']{\overline{Y}})$ is not dense in  $\CTt{\overline{Y}_1,P'}$).
Recalling that irreducible components of  $\CT{\overline{Y}_1}\times_{\tm^*/W_1} \tm^*$ have the same dimension and 
 each one contain $L_1\CN[\overline{P}_u]{\overline{Y}}$ as a dense subset for some $Y$, we see that there is only one such component.
\end{proof}

\begin{lemma}\label{MAIN_comp} Let $Y\in \fB_0(X)$. Then the equality of the closures of $G\NT(Y,P')$  and $G\NT(X,P')$ which are components of $\CTt{X,P'}$
  implies $Y=X$.
\end{lemma}
\begin{proof} If the equality holds then for  general $\xi\in \N(X,P')$ and $\eta\in \N(Y,P')$ we have $g(\xi,\Mom(\xi)|_{\s})=(g\xi,\Mom(\xi)|_{\s})=(\eta,\Mom(\eta)|_{\s})$,
By multiplying $\xi$ and $\eta$ by the elements of $P$
we can assume that  $\Mom(\xi)\in\ax^{pr}$ and $\Mom(\eta)\in \tm$. The above equalities give
 $\Mom(g\xi)=\Mom(\eta)=\Mom(\xi)\in \ax^{pr}$, i.e. $g$ is the stabilizer of the element of $\ax^{pr}$ that implies
 $g\in \LX$. Since $\LX$ normalizes a general $P$-orbit in $X$ we have $(g\xi,\Mom(\xi)|_{\s})\in\NT(X,P')$.
\end{proof}

\begin{lemma}\label{MAIN_comp_U} Let $Y\in \fB_0(X)$. Then the equality  of the closures of  $G\NT(Y,P_u)$ and $G\NT(X,P_u)$  which are components of $\CTt{X,P}$  implies $Y=X$.
\end{lemma}
\begin{proof} If the equality holds then for  general $\xi\in \N(X)$ and $\eta\in \N(Y)$ we have $g(\xi,\Mom(\xi)|_{\lv})=(g\xi,\Mom(\xi)|_{\lv})=(\eta,\Mom(\eta)|_{\lv})$,
	Again after multiplication by the elements of $P$ we can assume that  $\Mom(\eta)\in \tm$ and $\Mom(\xi)\in\tm^\bot$. Thus
	$\Mom(g\xi)=\Mom(\eta)=\Mom(\xi)\in \lv$, that implies
	$g$ fixes $\Mom(\xi)\in \tm^\bot$ and since $\xi$ is generic $g\in Z_G(\tm^\bot)$. Since $\tm^\bot \supset \ax^{pr}$ then $g\in \LX=Z_G(\ax^{pr})$.   Since $M$ normalizes a general  $P_u$-orbit in $X$ we have $(g\xi,\Mom(\xi)|_{\lv})\in\NT(X,P_u)$.	
\end{proof}

\begin{lemma}\label{sU_action} Let $Y\in \fB_0(X)$ such that $Y_1:=P_1Y\neq Y$, then  the raise from  $Y$ to $P_1Y$ is of type $(U)$. 
The action of $s\in W_{P_1}$ on the set of  components of $\CTt{X,P}$ (resp. $\CTt{X,P'}$) permutes the components corresponding to 
	$Y$ and $Y_1$.
\end{lemma} 
\begin{proof} The varieties  $\CN[P_u]{Y}$ and $\CN[P_u]{Y_1}$ (resp. $\CN[P']{Y}$ and $\CN[P']{Y_1}$) provide
the  components of $\CTt{X,P}$ (resp. components of $\CTt{X,P'}$).
By the reduction to $Y_1/(P_1)_u$ described in Lemma \ref{not_dense} and by Lemma \ref{MAIN_comp} and \ref{MAIN_comp_U} these components are distinct.
Moreover since $Y_1/(P_1)_u$ is the variety for the group $L_1$ of split rank $1$, there are no other subsets of  $Y_1/(P_1)_u)$ of maximal $k$-rank and $k$-complexity.
Since  $s$ acts nontrivially on  $\overline{\s}$, it must permute these components. 
\end{proof}

To prove $(ii)$ let us notice that since $P'$-orbits in $Y$ are normalized by $T_Y$ we have  $\Mom(\CN{Y,P'})\subset \aY{Y}+\p_u$. Since $G\CN{Y,P'}$ is dense in  $G\CN{X,P'}$
the general element  in $\Mom(\CN{Y,P'})$ is conjugate to $G\ab^{pr}$ and in particular semisimple. All semisimple elements of  $\aY{Y}+\p_u$ are $P$-conjugate to the elements of $\aY{Y}$
 This implies  ${\Mom(\CN{Y,P'})}\subset P\aY{Y}$.
The map $\Mom |_\s:\CN{Y,P'}\rightarrow \s^*$ factors through $\CN{Y,P'}\rightarrow \CT{Y} \rightarrow \s^*$, where the first map is a quotient by 
sub-bundle $\CNB{X}{Y}\subset \CN{Y,P'}$.
Moreover by Lemma \ref{restr_conorm} for a $\chi \in\Ch(Y)$ we have $f_\chi \in \kkk(Y)_\chi^{(P)}$, thus $df_\chi\in \CT{Y,P'} $ provides a rational section of the quotient of $\CN{Y,P'}$ by $\CNB{X}{Y}$. The points  of this section map
into $\chi+\p_u$ under $\Mom$. From the above we get $\overline{\Mom(\CN{Y,P'})}= P\aY{Y}$.

Let us notice that $P_u$-orbits in $Y$ are normalized by $T_Y$ and by $P_u\LXoo$. In particular  $\Mom(\CN{Y,P_u})\subset \aY{Y}+\Ad(w)\lve+\p_u$, where $\aY{Y}=w\ab$
and $w\lve\subset \lv$ (since $w\in N_G(L)$). This family has the section $w\K$ which is $wLw^{-1}$-invariant, has stabilizer of general position $L(Y)=wL(X)w^{-1}$ and the slice $wZ$. 
 Since $G\CN{Y,P_u}$ is dense in $\CT{X}$   the general element $\xi\in \Mom(\CN{Y,P_u})$ is semisimple.
 By Lemma \ref{adapt} the  element $\xi$  is $P_u$-conjugate to  its projection to $\aY{Y}+w\lv_{\eff}$. Moreover 
  the image of the element from  $\lv$  and its  $P_u$-conjugate have the same projection to $\lv$.
 This implies that image of moment map of $\CN{Y,P_u})$ as well as its projection to  $\lv$ can be recovered from the moment map of
 the restriction of $\CN{Y,P_u}$ to $w\K$. The map $\Mom[\lv]$  factors as $\CN{Y,P_u}|_{w\K}\rightarrow T^*(w\K)\rightarrow \lv^*$ 
 whose image is $w\Mom[\lv](\CT{X,\el})$ gives the equality  $\overline{\Mom(\CN{Y,P_u})}=Pw\Mom[\lv](\CT{X,\el})$ and proves $(ii)$.
\end{proof}

\begin{corollary} Let  the raise from $Y_1\in \fB_0(X)$  to $Y=P_1Y_1$ be not of type $(U)$. Then the transform $s_1w\K_\xi$ of the 
section $w\K_\xi$ for sufficiently general $\xi$ is also a section for $Y$.
\end{corollary} 
 \begin{proof} 
 By Lemma \ref{not_dense} $G\CNt{Y,P'}$ is dense in some component of $\CTt{X,P'}$   and $G\CNt{Y_1,P'}$ is not dense in any component
 for all $Y'\in \fB_0(Y)$. Since $\Mom(\CN{Y,P'})=P\ab_Y$, the the same holds for the $G$-span of the 
   lift $\widetilde{C}_Y\subset \CTt{X,P'}$ of  $C_Y:=\CN{Y,P'}\cap \Mom^{-1}(\ab_Y)$.
Since $s_1C_Y$ is contained in the union of $\CN{Y',P'}$ for $Y'\in \fB_0(Y)$, and in the same time $G\widetilde{C}_Y$ is dense in $\CTt{X,P'}$,
the latter implies that $s_1C_Y=C_Y$ as well as $s_1\ab_Y=\ab_Y$. 
Since $C_Y=\bigcup_{\xi\in \ax^{pr}}w\K_{\xi}$, we have $s_1w\K_{\xi}=w\K_{w^{-1}s_1w\xi}$. 
 \end{proof}

\begin{proof}[Proof of Proposition \ref{Wprop}] 
Given $Y\in\PS(X)$ we associate the corresponding $k$-irreducible component by  taking the $G$-span of the
  lifting of $\CN{Y,P_u}$ (resp. $\CN{Y,P'}$) to $\CTt{X,P}$ (resp. $\CTt{X,P'}$).
By Lemma  \ref{sU_action}  and \ref{not_dense}  the action of generators of $F(W_k)$  on $Y\in\PS(X)$ is consistent with the action
of generators of $W_k$ on the set of corresponding components, that gives homomorphism.  By Lemma   \ref{MAIN_comp} (\ref{MAIN_comp_U}) the map from the set of principal families $\PS(X)$
to the components of $\CTt{X,P}$ (resp. $\CTt{X,P'}$)  is injective
and by Proposition \ref{dim_comp} it is surjective. 
\end{proof}

\section{Functoriality properies of the Weyl group actions for field extensions.}\label{s:functoriality}

In this part we are interested in studying functoriality properties of the Weyl group actions under the field extensions $[E:F]$. 
All the objects defined before
will have additional index, for example, $P_E$ denotes a minimal parabolic defined over $E$ and $P_E(E)$ denotes  the set of its $E$-points.
The set of  principal closed $P_E$-subvarieties 
 in $X$ is denoted by    $\PS(X,E)$, and the action of $w\in W_E$ on $Y\in \PS(X,E)$ is denoted by $w\circ_E Y$. For the subgroup $H\subset G$ defined over $E$ and normalized by $\AF{E}$ by $\Delta^{\pm}_{H,E},\Pi_{H,E}$ we denote the set of
positive, negative and simple $\AF{E}$-roots where the positivity is related to the choice of minimal parabolic $P_E$.
The functorial properties for the field extensions can be summarized in the following diagram.

\xymatrix{
                                                                                     & \CT{X}\times_{\tm^*/W}\s_E^*   \ar[dl] \ar[d]                  & \ar[d] \ar[l] \CT{X}\times_{\tm^*/W}\s_F^*                 &    \\
\CT{X}\times_{\tm^*/W}\tm^*   \ar[r]^{/W_{L_E}}  & \ar[d]^{/W_E}  \CT{X}\times_{\tm^*/W}\lv_E^*/\!\!/L_E \ar[r]    & \ar[d] ^{/W_F} \CT{X}\times_{\tm^*/W}\lv_F^*/\!\!/ L_F     & \\
                                                                                    &  \CT{X}\times_{\tm^*/W}\lv_E^*/\!\!/N_G(L_E) \ar[r] \ar[rd]                 & \CT{X}\times_{\tm^*/W}\lv_F^*/\!\!/ N_G(L_F)   \ar[d]       & \\
												&														&		\CT{X}												&
}

We shall define the action of some extension of $\widetilde{W}_F\subset W_E$ of $W_F$ on the set $\PS(X,E)$ by choosing a lifting of the generators of $W_F$ to $N_G(\AF{E})$ defined over $E$ and generating 
a subgroup $\widetilde{W}_F\subset W_E$ by these generators.
 This can be done as follows. 

Let  $s\in W_F$ be the generator.
  Over $E$ we may choose the minimal parabolic subgroup $P_E\subset P_F$ as the preimage  of a minimal parabolic subgroup in $L_F:=P_F/(P_F)_u$.
 Since $s$ normalize $L_F$, the group $sP_Es^{-1}\cap L_F$ is also a minimal parabolic $E$-subgroup of $L_F$ over $E$  and it is conjugate to $P_E\cap L_F$  by some $l\in L_F(E)$.  
 The maximal $E$-split tori in $P_E$ are conjugate by element of $P_E(E)$, so we may specify the choice of $l$ in such a way that $\widetilde{s}=ls$ normalizes $A_E$.
 This  defines $l$ up to the element  of $L_E=P_E\cap N_G(\AF{E})$ so $l$ defines a unique element of $W_E$. 
Since the lifting $\widetilde{s}$ normalize $P_E\cap L_F$ it also permutes the supminimal parabolic subgroups of $L_F$ containing  $P_E\cap L_F$.
These supminimal parabolics correspond to the  set of simple $\AF{E}$-roots and, in particular, 
 if $\alpha$ is the simple $\AF{E}$-root of $L_F$,  then $\beta:=s\alpha s^{-1}$ is also a simple $\AF{E}$-root. 

This can be restated as follows. Consider a simple root $\gamma \in \Delta_F$ and  let $(\AF{F})_\gamma$ be the torus 
that stabilizes the corresponding root subspace $\g_\gamma$. Considering the group $Z_G((\AF{F})_\gamma)$ instead of $G$ we can assume that $G$ is of semisimple split rank $1$.
Since $P_F$ and $P_F^-$ are opposite minimal parabolic subgroups defined over $F$ then they are conjugate by the element $s\in G$ that preserves $L_F$. We can assume that $s$
maps $(P_F)_u$ into $(P_F^-)_u$ and  $P_E\cap L_F$ maps into itself. 
As above, this defines the element of $\widetilde{s}\in W_E$ whose  the length in $W_E$ that is the number of positive 
$\AF{E}$-roots made negative by applying $s$ and this number is equal to the
number of $\AF{E}$-root subspaces in $(\p_F)_u$. By this definition we have $\widetilde{s}=w_{L_F,E}w_E$, where $w_E$ is the longest element of $W_E$
i.e. a unique element such that $w_E\Delta^+_E=\Delta^-_E$ and $w_{L_F,E}$ is the longest element for $L_F$ (related with the root system of $\AF{E}$) i.e.
 a unique element such that $w_{L_F,E}\Delta^+_{L_F,E}=\Delta^-_{L_F,E}$).


\begin{theorem}\label{E_F_compl}  The set $Y\in \PS(X,F)$ is $P_E$-invariant and provides the principal family 
of $P_E$-orbits, i.e. the element of $\PS(X,E)$.
\end{theorem}
\begin{proof} Let  recall that by Lemma \ref{stab} there is a sequence of supminimal parabolic subgroups $P_i$, defined over $F$  
 and a sequence of subvarieties $Y_i:=P_{F,i}\ldots P_{F,1}Y$ such that the raise from $Y_{i-1}$ to $Y_i$ is of type $(U)$ and $Y_l=X$.
Let us notice that if we take the  decomposition of minimal length $\widetilde{s}_i=s_{\alpha_{l_i}}\ldots s_{\alpha_1}$  for the lifting of $s_i$ to $W_E$, where $\alpha_j\in \Pi_E$ are simple $\AF{E}$-roots,
the raise from $Y_{i-1,j}=P_{E,\alpha_j}\ldots P_{E,\alpha_1}Y_{i-1}$ to  $Y_{i-1,j+1}=P_{E,\alpha_{j+1}}\ldots P_{E,\alpha_1}Y_{i-1}$ is also of type $(U)$. The raise of type $(U)$
for $P_F$-invariant subvarieties  is characterized by the property 
there exists a generic point $y\in Y_i$ such that the image of $(P_F)_y$ in  contains the unipotent subgroup $(U_F)_{\gamma_i}$ corresponding to the root $\gamma_i\in \Pi_F$. In the same time  the words 
$w_j:=s_{\alpha_j}\ldots s_{\alpha_1}$ have minimal length considered  as the elements of $W_E$. The roots  $\alpha_i$ are the roots of the semisimple part of $P_{E,i}$. 
Since $\ell(s_{\alpha_{j+1}}w_{j})=\ell(w_j)+1$ the $w_j(\Delta^+_E \cap w^{-1}\Delta^-_E)$ consists of a  root $\alpha_{j+1}$ which is a unique $A_E$-root that changes the sign after application of $s_{\alpha_{j+1}}$. Since the set of $A_E$-roots from the set $\Delta^+_E \cap w^{-1}\Delta^-_E$ are precisely the $A_E$-roots of the group $(U_F)_{\gamma_i}$,
 the stabilizer of $w_jy$ also contains the  root subspace $(U_F)_{\alpha_{j+1}}$ and the corresponding raise from $Y_{i-1,j}$ to $Y_{j+1,i-1}$ is of type $(U)$.
Thus we get that $Y$ is connected to $X$ by the sequence of raises defined over $E$ of type $(U)$.
Since the $(P_E)_u$-complexity does not increase after the raise of type $(U)$  this proves that  $Y$ also has maximal $(P_E)_u$-complexity i.e. $Y\in \PS(X,E)$.
\end{proof}

Let us fix the following notation till the end of the section. Let $\tilde{s}\in W_E$ be the lift of simple reflection $s\in W_F$ corresponding to the 
simple root $\gamma \in \Pi_F$. Let $P_{F,\gamma}$ be a corresponding supminimal parabolic. Let us take $\alpha \in \Pi_{L_F,E}$  and
$\beta := \tilde{s}\alpha \tilde{s}^{-1}\in \Pi_{L_F,E}$.

\begin{proposition} Let $Y_F\in \PS(X,F)$, $Y\in \PS(Y_F,E)$ and $P_{F,\gamma}$ be a supminimal parabolic such that the raise
from $Y_F$  to $P_{F,\gamma}Y_F$ is of type $(U)$. Then the word $\widetilde{s}_\gamma=s_{\alpha_{l}}\ldots s_{\alpha_1}$ defines the 
decomposition $P_{F,\gamma}=P_{F,\alpha_{l}}\ldots P_{F,\alpha_1}$ that allows to describe
 the raise from $Y$ to $P_{F,\gamma}Y$ as the composition of the raises  of type $(U)$.
\end{proposition} 
\begin{proof}
Let us notice that by 
 \cite[Prop. 4.13]{KK} we have  $c(Y/(P_E)_u)\leq c(P_FY/(P_E)_u)$ and the equality can take place
only if $Y_F=P_FY$. Let us notice that the raise from $Y_F$ to $P_{F,\gamma}Y_F$ is of type $(U)$, iff the stabilizer of some $y\in Y_F$ contains  
the unipotent subgroup $(U_F)_{\gamma}$. This property is preserved after conjugation by the element of $l\in L_F$, and since $Y_F=L_FY$
this property holds for some point of $Y$ as well.  Now we can finish the proof similar to the proof of Theorem \ref{E_F_compl}.
\end{proof}



\begin{proposition} Let $Y_F\in \PS(X,F)$ and $Y\in \PS(Y_F,E)$. Let $\tilde{s}\in W_E$ be the lift of simple reflection $s\in W_F$ corresponding to the 
simple root $\gamma \in \Pi_F$, $\alpha \in \Pi_{L_F,E}$  and
$\beta = \tilde{s}\alpha \tilde{s}^{-1}$.
   Then the type of the raise  from $Y$ to $P_{E,\alpha}Y$ is the same as the type of the raise $\tilde{s}\circ Y$ to $P_{E,\beta}(\tilde{s}\circ Y)$ 
\end{proposition}
\begin{proof} Given the sets $Y_F$ and $P_{F,\gamma}Y_F \in \PS(X,F)$ we can pass to the some open subset of $Y_F$, in order to assume existence of  the $L_F$-varieties $Y_F/(P_F)_u$ and $P_{F,\gamma}Y_F/(P_F)_u$. Let us note  that for $y\in Y_F$ and its image $\overline{y}\in Y_F/(P_F)_u$, the  stabilizer  $(L_F)_{\overline{y}}$ is equal to the image of $(P_F)_y$ in  $L_F$. If we fix the section  $\K_E$ defined over $E$ for the subset $Y\in \PS(Y_F,E)$,  we get the following diagram:

\xymatrix{
&	\ar[r] P_{F,\gamma} Y		&\ar[r]	P_{F,\gamma}P_{E,\alpha} Y=P_{E,\beta}P_{F,\gamma} Y	& P_{F,\gamma}Y_F/(P_F)_u \\
&	\ar[r]\ar[u] Y		&\ar[r] \ar[u] P_{E,\alpha} Y		& \ar[u] Y_F/(P_F)_u
}

Which can be  rewritten in terms of sections: 

\xymatrix{
&	\ar[r] 	P_F\widetilde{s}\K_E	&\ar[r]	P_F\widetilde{s}s_{\alpha}\K_E=P_Fs_{\beta}\widetilde{s}\K_E	& P_{F,\gamma}Y_F/(P_F)_u \\
&	\ar[r]\ar[u] P_F\K_E		&\ar[r] \ar[u] P_Fs_{\alpha}\K_E		& \ar[u] Y_F/(P_F)_u
}

The type of the raise from $Y$ to in $P_{E,\alpha} Y$ is the same as the type of the raise for their images in $Y_F/(P_F)_u$
and is defined by the type of the stabilizer $(L_F)_{\overline{y}}\cap P_{E,\alpha}$ (where $y\in \K_E$ and $\overline{y}$ is its image in $Y_F/(P_F)_u$).
in the same time the type the raise from  $P_{F,\gamma} Y$	to $P_{E,\beta}P_{F,\gamma} Y$
is defined by the type of the stabilizer $(L_F)_{\widetilde{s}\overline{y}}\cap P_{E,\beta}$ which is obtained from $(L_F)_{\overline{y}}\cap P_{E,\alpha}$
via conjugation by $\widetilde{s}$.  
\end{proof}

\begin{remark} Assume that for the element $\widetilde{s}\in W_E$ (which is a lift of $s\in W_F$) and $Y_F\in \PS(X,F)$ we have $\widetilde{s}\circ Y_F\subset Y_F$.
Then $\widetilde{s}$  acts on the graph of $P_E\cap L_F$-orbit inclusions for the $L_F$-variety $Y/(P_F)_u$ by the automorphism
 that $\widetilde{s}$ induces on $L_F$ (on the set of simple roots of $L_F$) and  $\widetilde{s}\circ Y_F= Y_F$.
\end{remark}

The set  $Y_F$ is irreducible over $E$ since it contains  a Zariski
 dense subset of $F$-points. In particular it 
contains a unique subset of $(P_E)_u$-orbits of maximal dimension that gives the corresponding
subset of $(P_E)_u\cap L_F$-orbits in $Y_F/(P_E)_u$ which maps into itself by automorphism $\widetilde{s}$. So we have the following corollary:

\begin{corollary} Suppose that the action of $s\in W_F$ on the irreducible over $E$ subset $Y_F\in \PS(X,F)$ is not of type $(U)$. Then the $W_E$-action of the lift $\widetilde{s}\in W_E$  
 on the  subset $Y_F$ considered  as the element of $\PS(X,E)$ fixes $Y_F$.
\end{corollary}

Combining the above results we get the following theorem:

\begin{theorem}  The action of $W_E$ on $\PS(X,E)$ induces the action of $W_F$ on $\PS(X,F)$ via the lift defined above.
\end{theorem}

\section{Generators of the little Weyl group.}\label{s:generators}

For  $Y\in \PS(X)$ let $\WE{Y} \subset W$ be the stabilizer of $Y$ with respect to the Weyl group action on the principal $P$-invariant closed irreducible subsets of $X$.
The proof of the following theorem is the smart combination of the proofs  of  M.Brion \cite[Prop.4]{BBORB} and F.Knop  \cite[$\S$7]{BORB} with some new ingredients.

\begin{theorem} The subgroup $\WE{X}\subset W_k$ is generated by the reflections $s_\alpha\in W_k$ or the products of reflections $s_\alpha s_\beta$ such that $s_\alpha,s_\beta\in W_k$,
 $\alpha \ann \beta$ and $\alpha+\beta\in \ab$ is not a root of $\Delta_k$.
\end{theorem}
\begin{proof}
We shall first make several observations that are essential for the proof of the theorem.
Let $Y\in  \PS(X)$ and $\alpha, \beta \in \Pi$ be the simple  roots such that $\dim P_\alpha Y> \dim Y$ and $\dim P_\beta Y> \dim Y$.
Then since the $k$-rank, $k$-complexity and $k$-homogeneity are maximal, then the raise from $Y$ to $P_\alpha Y$ (resp. to $P_\beta Y$) is of type $(U)$, that also implies that for $P_\alpha$ 
 (resp. for $P_\beta$) we have
$$\dim P_\alpha Y=\dim Y+\dim P_\alpha/P.$$ 
 Let $P_{\alpha\beta}$ be the subgroup
of $G$ generated by $P_\alpha$ and $P_\beta$.  Consider the $P$-invariant subvariety $Y_{m}^\alpha:=\underbrace{P_\alpha P_\beta\ldots}_{m-terms}Y$ of $X$ obtained
 by $m$-subsequent applications of $P_\alpha$ and $P_\beta$. Let us notice that if $Y_\alpha$ is  also $P_\beta$ stable  then $Y_{m}= P_{\alpha\beta}Y$.
Similarly define $Y_{m}^\beta:=\underbrace{  P_\beta P_\alpha \ldots}_{m-terms}Y$
 
 For a simple $k$-root $\alpha$ let us denote $n_\alpha=\dim P_\alpha/P$.
 This implies that unless $Y_{m-1}^\alpha=P_{\alpha\beta}Y$ we have 
 $\dim Y_{m}^\alpha=\dim Y_{m-1}^\alpha+n_\alpha$ and 
$\dim Y_{m}^\beta=\dim Y_{m-1}^\beta+n_\beta.$
 Let us notice that there exist $m_\alpha,m_\beta$ for which:
 $$
 P_{\alpha\beta}Y=\underbrace{P_\alpha P_\beta\ldots }_{m_\alpha-terms}Y=\underbrace{P_\beta P_\alpha\ldots}_{m_\beta-terms}Y
 $$

Comparing the dimensions of the terms of this equality, we get that $m_\alpha=m_\beta$, and this number should be even in the case $n_\alpha\neq n_\beta$.

This equality implies that $(s_\alpha s_\beta)^m\cdot Y=Y$, and in particular $$(s_\alpha s_\beta)^m\Ch(Y)=\Ch(Y).$$ 
Assume that $(s_\alpha s_\beta)^m\neq e$.
Since $P_\alpha$ do not normalize $Y$ and the raise is of type $(U)$ by Claim \ref{degen_weight}  we have $\alpha\notin \DY{Y}$.
Thus  the roots $\alpha,\beta$ are not orthogonal to $\aY{Y}$.

Let $\Rad(P_{\alpha\beta})$ be the radical of $P_{\alpha\beta}$,
$W_{\alpha\beta}\subset W$ be the  subgroup generated by $s_\alpha$ and $s_\beta$,
and $G_{\alpha\beta}:=P_{\alpha\beta}/\Rad(P_{\alpha\beta})$.

\begin{lemma} There is an isomorphism of the set $\PS(P_{\alpha\beta}Y,P)$ and the set
 $$\PS(P_{\alpha\beta}Y/\Rad(P_{\alpha\beta}),P/\Rad(P_{\alpha\beta}))$$ equipped with the action of
$W_{\alpha\beta}$.
 \end{lemma}

This lemma reduces our situation to the study of  homogeneous   space $X:=G_{\alpha\beta}/H$. For brevity put $G:=G_{\alpha\beta}$.

\begin{lemma} The plane $H_{\alpha\beta}$ is not contained in $\aY{Y}$ unless  $G/H$ is horospherical and $(s_\alpha s_\beta)^m=e$.
\end{lemma}
\begin{proof} After taking quotient by the  radical of  $P_{\alpha\beta}$ we reduce the situation to the homogeneous
space of the group of rank $2$ that we denote by $G$. 
Since $Y$ is not stable under $P_\alpha,P_\beta$ there are no orbits that are raised to $Y$. In particular, all $P$-orbits of the points from $G/H(k)$ have maximal rank equal to $2$, and in particular the 
stabilizer $gHg^{-1}\cap P$ does not contain nontrivial $k$-split torus for all $g\in G_k$. Since all $k$-split torus is conjugate to the subtorus of $P$ by some element $g\in G_k$. 
This implies that $H$ do not contain $k$-split tori and in particular $H/\Rad_u H$ is anisotropic.  
Moreover since the raise from $Y$ to $P_\alpha Y$ is of type $(U)$,  $H$ has a subgroup isomorphic to $\Bbb G_a$, that implies $\Rad_u H$ is non-trivial.
 We can assume
$H$ belong to one of the parabolic groups $P_i:=P,P_\alpha,P_\beta$ and $H_u$ belong to the corresponding unipotent radical $(P_i)_u$. Also we can assume after conjugation by 
$p\in P_i(k)$ that a Levi subgroup $\Lv_H$ of $H$ (which is elementary) is contained in the Levi subgroup $\Lv_i$ of $P_i$.  
  Consider the equivariant map $G/H\rightarrow G/P_i$. The $\PU$-orbits on $G/H$ 
 that map to  $\PU wP_i/H$ (where $w\in W_k$) correspond to orbits  $w^{-1}\PU w\cap \PU$ on  $\PU/H_u\times \Lv_i/\Lv_H$
  (which correspond to orbits of stabilizer $w\PU w^{-1}\cap \PU$ on fiber  $wP_i/H$).
  For the open cell $ \PU w_0P_i/H$ this stabilizer is trivial and the dimension of such family is equal to $\dim(P_i/H)$.  
  Since the raise from   $Y_{m-1}^\alpha$ and from $Y_{m-1}^\beta$ to $G/H$ is of type $(U)$ then the dimension of general $\PU$-orbits in  $Y_{m-1}^\alpha$ and $Y_{m-1}^\beta$ is strictly smaller than $\dim \PU$ so these divisors do not intersect $\PU w_0P_i/H$, and moreover $\codim_X Y_{m-1}^\alpha=\codim_{G/P}Bs_\beta w_0P/P=n_\beta$.  
If  $m=1$, then $n_\alpha=n_\beta$ and $Y$ is equal to preimage of one of the Schubert variety of codimension $n_\alpha$
 which is impossible since $Y$ is neither $P_\alpha$ nor  $P_\beta$-stable and the preimage of the Schubert divisor is stable with respect to one of these groups.
 
 In the case  $Y_{m-1}^\alpha$ and $Y_{m-1}^\beta$ are $k$-dense and have codimension $n_\beta$ and $n_\alpha$ in $Y$ respectively and do not contain in the closure of each other.
  On the other hand their images to $G/P_i$  coincide with  $k$-dense Schubert varieties  $\PU wP_i/H$ for $w\in W_k, w\neq e$.
  But the minimal codimension   for such Schubert varieties that do not lie in the closure of larger variety with this property
    is equal to $n_\alpha$ and  $n_\beta$ which is also equal to codimension in $X$
  of their preimages  under the equidimensional map to $G/P_i$.
    
  Thus  $P_i=P$ and $Y_{m-1}^\alpha$ and $Y_{m-1}^\beta$  are equal to  
  $\PU s_\beta w_0P/H$ and $\PU s_\alpha w_0P/H$ which are the preimages of different $k$-dense Schubert varieties.  
 Since the raise of type $(U)$ preserve the $\PU $-complexity i.e. $c(X/\PU)=c(Y_{m-1}^\alpha/\PU)=c(Y_{m-1}^\beta/\PU)$ then the complexity of the actions  
of  $w_0^{-1}s_\alpha \PU s_\alpha w_0\cap \PU$ and  $w_0^{-1}s_\beta \PU s_\beta w_0\cap \PU$ (which are the root subgroups ${(\PU)}_\beta$ and ${(\PU)}_\alpha$)  on $P/H$ 
is equal to $\dim (P/H)$, that implies the  triviality of these actions. 
 
  In particular the  triviality of the action of root subgroups ${\PU}_\beta$ and ${\PU}_\alpha$ implies the triviality of the action of $\PU$ on $P/H$, thus $\PU=H_u$.
 In this case $\W{(X)}$ is trivial  since the families of orbits of maximal  $\PU$-complexity equal to $\dim(\Lv/\Lv_H)$ 
 form the subsets to $\PU wP/H$ for $w\in W_k$, and $W_k$ acts transitively on them 
(by \ref{MAX_RC}). Thus we have $(s_\alpha s_\beta)^m=e$ that contradicts with our assumptions.
\end{proof}

The element $(s_\alpha s_\beta)^m$ is the rotation in the two dimensional plane $H_{\alpha\beta}$  spanned by $\alpha$ and $\beta$ 
which preserves  $\aY{Y}$ (and fixes $H_{\alpha\beta}^\bot$). By the previous Lemma and since $\alpha,\beta$ are not orthogonal to $\aY{Y}$ the 
subspace $\aY{Y}\cap H_{\alpha\beta}$ is one dimensional  and fixed by rotation $(s_\alpha s_\beta)^m$ that implies:

\begin{proposition}  $(s_\alpha s_\beta)^m=-\Id|_{H_{\alpha\beta}}$ and the restriction of $(s_\alpha s_\beta)^m$ to $\aY{Y}$ is the reflection in the hyperplane
orthogonal to $H_{\alpha\beta}\cap \aY{Y}$.
\end{proposition}

This  can happen only in the following cases:
$(i)$ $m=1$, $\alpha\ann\beta$,
$(ii)$  $m=2$, $\alpha, \beta$ generate the root system of type $B_2$, $BC_2$,
$(iii)$ $m=3$, $\alpha, \beta$  generate  the root system of type $G_2$.

\begin{proposition} 
Let $G$ be a group of rank two,  then we have the following possibilities:
either $G/H\cong PGL_2\times PGL_2/PGL_2$ and $\ab_{X}=\langle \alpha+\beta \rangle$,
or   $m>1$ and  $(s_\alpha s_\beta)^m$ is equal to a  reflection  $s_\gamma \in W$, where $\gamma \in \Xi(Y_i)$.
\end{proposition}
\begin{proof} 

Let us prove first that either $G/H$ satisfies the conditions of proposition or it  is spherical.
  The subset $G(k)x_0\cong G(k)/H(k)$ for $x_0\in X(k)$ is Zariski
 dense in $G/H$ and correspond to the set $G(k)/P(k)$
see Lemma \ref{P_G_H}. Since a general $P$-orbit has $k$-rank $r$ there exists
 a  split torus $A_0$ of dimension $\rk_k(G/H)-r$ such that for each $P$-orbit of  $x=gx_0$, after conjugation by $p\in P(k)$ we have $A_0\subset gHg^{-1}\cap P$.
We can also assume that $A_0$ is contained in some fixed maximal split torus $A_H$ of $H$.
Let us notice that the set of $g\in G$ such that $g^{-1}A_0g\subset A_H$ is closed and in particular has finite number of connected
	components. Since the a subtorus of $A_H$ cannot vary continuously, we can assume that there exists $A'\subset A_H$ which is $H(k)$-conjugate
	to $g^{-1}A_0g$ where $g\in G(k)$. In this way for a generic  $H$-orbit we have associated a subset $Y'\in (G/P(k))^{A'}$, and for those $H$-orbits
	that  contain a point from $G/P(k)$ we have associated a point in $Y'(k)\in(G/P(k))^{A'}$. In particular if $G/H$ is not spherical $Y'(k)$ is infinite.

Assume that $\Ch(X)$ is not generated by a root. Then  $(G/P(k))^{A'}=\{sP,eP\}$  since this set of points coincide with the fixed points of $A$.
In this case  $G/H$ is spherical.

 In the case,  assume that $[\LXk,\LXk]$ be nonanisotropic. Then the root system $\Delta_{\LXk}$ is generated by $\alpha$ or $\beta$,  $\ab_{X_{\alpha\beta}}$ 
 is generated by the corresponding orthogonal root
and $(s_\alpha s_\beta)^m$ is the reflection with respect to this root. 

When $[\LXk,\LXk]$ is anisotropic, by the local structure theorem  $P_u$ acts freely on the points of open $P$-orbit. The following lemma finishes the proof of the proposition. 

\begin{lemma} Let $G_{\alpha\beta}/H$ be a  spherical  homogeneous space   with  a locally free action of 
 $\PU$ and such that $H$ contains regular one-dimensional split torus $A_0$
 (i.e. $Z_{G_{\alpha\beta}}(A_0)=\Lv$).
Then either $G_{\alpha\beta}$ is of type $A_1+A_1$ and $\ab_{X_{\alpha\beta}}=\langle \alpha+\beta \rangle$  or $G_{\alpha\beta}/H$ is horospherical.  
\end{lemma}
\begin{proof}
By the local structure theorem we have a choice of $P$ and $H$ such that $\p+\h=\g$ and $\p_u\cap \h=\{0\}$. 
Moreover we may assume that  $\Xi({X_{\alpha\beta}})$ is contained in the positive Weyl chamber. 
This can be done by passing to quasiaffine  $G/H$ via taking an affine  cone over $G/H$ and extending
$G_{\alpha\beta}$ and $H$ by dilatations.
The torus $\ab_0$ is regular and orthogonal subspace $\Xi({X})$ belong to the interior of the Weyl chamber. Denote by $\lambda_0$ a one-parameter subgroup corresponding to $\tm_0$
In the case when $A_0$ has different weights on all the roots subspaces of $\g_{\alpha\beta}$,  $\h$ is normalized by $A_0$ and it is the direct sum of $A_0$ subspace of $\lv$ and subspaces 
of distinct root subspaces. The the above equalities imply $\h=\lv_{x}+\ab_0+\p_u^-$ for $\lv_x\in [\lv,\lv]$. 
The direct check shows that $\lambda_0$ has the same pairing  with different roots only in the following cases: 
$(i)$ $\Delta$ of type $A_1+A_1$, $\ab_{X_{\alpha\beta}}=\alpha+\beta$, and $\lambda_0=\alpha-\beta$
has the same pairing with the roots $\alpha$ and $-\beta$. This is the precisely the statement of the Lemma. 
$(ii)$ $\Delta$ of type $B_2$, $\ab_{X_{\alpha\beta}}=2\beta+3\alpha$, and $\lambda_0=\alpha-\beta$
has the same pairing with the roots $-\alpha-\beta$ and $2\alpha+\beta$ and the same holds for the opposite pair, $(iii)$ $\Delta$ of type $G_2$, $\ab_{X_{\alpha\beta}}=5\alpha+3\beta$, and $\lambda_0=\alpha-\beta$
has the same pairing on the roots $-2\alpha-\beta$ and $3\alpha+2\beta$.
 But in the cases $(ii)$ and $(iii)$ the pairing of $\lambda_0$ with the roots $-\alpha$, $-\beta$ are distinct and differ from 
the pairings with other roots.  Since $\p+\h=\g_{\alpha\beta}$,  $\h$ and $\p$ is normalized by $T_0$ and the  root subspaces  of $-\alpha$, $-\beta$  do not belong to
$\p$ they are contained in $\h$. Since subspaces  $\g_{-\alpha}$ and $\g_{-\beta}$ generate $\p_u^-$ this proves the lemma. 
\end{proof}

Since  the space $G/H$ is horospherical,  it has a trivial little Weyl group and also trivial $W_{(X)}$ due to triviality of $[L,L]$. Thus $(s_\alpha s_\beta)^m=1$, which finishes the proof of proposition.
 
 \end{proof}

Let $w\in \WE{X}$ and let $s_{\alpha_\ell}\ldots s_{\alpha_1}$ be a reduced decomposition of $w$. Let
us put  $w_i=s_{\alpha_i}\ldots s_{\alpha_1}$ and $Y_i=w_i\cdot X$.

We finish the proof of the theorem by the induction on the following three parameters:
\begin{itemize}

\item[$\bullet$]   $\mcd=\{\text{the maximal codimension of} \ Y_i \ \text{in}\  X\ \text{for all}\ i\}$,
\item[$\bullet$]  $ \nd=\{\text{the number of} \ Y_i | \codim_{Y_i}   X=\mcd\}$,
\item[$\bullet$]   $\ld=\{ \text{length of  the given decomposition of} \ w\}$.

\end{itemize}

On each step of the induction we shall either decrease $\mcd$
(possibly increasing  $\ld$ and $\nd$),  or decrease $\nd$ not changing $\mcd$ (possibly increasing $\ld$), or  decrease $\ld$, not changing $\mcd$ and $\nd$.

Consider the minimal $i$ for which $Y_i$ has maximal codimension in $X$.

{\bf Case 1} Let $P_{\alpha_{i+1}}Y_i=Y_i$ (in particular $s_{\alpha_{i+1}}\in \WE{Y}$) then we have two cases:
there either exists or not exisits $Y'\in \B_0(X)$ with $\dim Y'=\dim Y_i-1$, such that $P_{\alpha_{i+1}}Y'=Y_i$.
In the second case $P_{\alpha_{i+1}}$ is the normalizer of the general $B$-orbits in $Y_i$ and $\alpha_{i+1}\in \DY{Y_i}$.
In the first case  the type of  the pair $(Y_{i},Y')$  is different from the type $(U)$, which implies that
the vector space spanned by the character lattice $\Ch(Y_i)$ is $s_{\alpha_{i+1}}$-stable.
Let us notice that $\alpha_{i+1}\notin  \DY{Y_i}$ so $\alpha_{i+1}$  is not orthogonal to $\Ch(Y_i)$ (due to non-degeneracy assumption).
In particular this implies that $\alpha_{i+1}\in \Ch(Y_i)$.
Since $w_i^{-1} s_{\alpha_{i+1}}w_i\in \WE{X}$,  the element
 $$w(w_i^{-1} s_{\alpha_{i+1}}w_i)=s_{\alpha_\ell}\ldots \widehat{s_{\alpha_{i+1}}}\ldots s_{\alpha_1}.$$
also belongs to $\WE{X}$ and has the smaller length than $w$.

{\bf Case 2}  Let  $\dim P_{\alpha_{i+1}}Y_i>\dim Y_i$. From the first part of our proof we see that there
exists $m$ such that $(s_{\alpha_{i+1}}s_{\alpha_{i}})^m\in \WE{Y_i}$
and either $s_\beta= (s_{\alpha_{i+1}}s_{\alpha_{i}})^m$ is the reflection of $W$
and $\beta \in \Ch(Y_i)$
 or $m=1$ and $\alpha_{i+1}$ and ${\alpha_{i}}$ are
orthogonal roots whose sum is not a root of $\Delta$. In the last case  $\alpha_{i+1}+{\alpha_{i}}\in \Ch(Y_i)$.
 As before $w_i^{-1} (s_{\alpha_{i+1}}s_{\alpha_{i}})^m w_i\in \WE{X}$, consider the element of $\WE{X}$
$$ww_i^{-1} (s_{\alpha_{i+1}}s_{\alpha_{i}})^m w_i= s_{\alpha_\ell}\ldots s_{\alpha_{i+2}} (s_{\alpha_{i}}s_{\alpha_{i+1}})^{m-1}s_{\alpha_{i-1}}\ldots s_{\alpha_1}.$$

 It will have the larger minimal $i$ for which $Y_i$ has maximal codimension in $X$ and the number of such $Y_i$ of maximal codimension will be decreased. This proves the theorem.
\end{proof}

\section{Monodromy action on the set of real orbits}\label{s:monodromy}

In this section we focus again on the case when $k=\Bbb R$, and the variety $X$ is homogeneous and admits the so-called $k$-wonderful compactification. For  the definition of the concept of  $k$-wonderful variety and the required details see  \cite{KK}. Also we shall need the  following theorem:

\begin{theorem}{\cite[Thm.13.7]{KK}}
  Let $X=G/H$ be a homogeneous $\RR$-spherical
  variety. Assume that $X$ is $\RR$-wonderful and let $X\into X^{\rm st}$ be
  its  wonderful embedding. Let $Y\subseteq X$ be the
  closed $G$-orbit isomorphic to $G/\PX^-$. Then

\begin{enumerate}

\item\label{it:Rwonderful2} $X\st(\RR)$ is a compact connected
  manifold.

\item\label{it:Rwonderful3} $Y(\RR)$ is the only closed
  $G(\RR)^0$-orbit of $X\st(\RR)$. In particular, it is connected and
  $G(\RR)$-stable.

\item\label{it:Rwonderful1} $P(\RR)^0$ has at most $2^{\rk_\RR X}$
  open orbits in $X(\RR)$ (or, equivalently, in $X\st(\RR)$). They all
  contain $Y(\RR)$ in their closure. 
\end{enumerate}
\end{theorem}

Let us also recall that we can calculate the topological fundamental group of $Y(\RR)$ using the fact that for the set $Y(\RR)$ we have a Bruhat  decomposition
which is also a cell decomposition of the corresponding $CW$-complex (cf. \cite{FUND} where $\pi_1(Y(\RR))$ is described by generators and relations). 
The generators of $\pi_1(Y(\RR))$ correspond to  one dimensional cells in $Y(\RR)$ which in turn correspond to
the simple roots $\alpha \in \Delta_k\setminus \Delta_{\LX}$ such that  $\dim P_\alpha/P=1$ (or equivalently $ \g_{2\alpha}=0$ and $\dim \g_\alpha=1$).
The loop  is represented by applying $\exp(\pi t Z_\alpha)$ to a $P$-fixed point $y_0\in Y$ for $t\in [0,1]$, where $Z_\alpha=e_\alpha+e_{-\alpha}$.

Let us notice that the set  $X(\RR)$ is obtained by throwing away
$\RR$-points of complete intersection divisor (its smooth components we denote by $E_i$) from the compact manifold $X\st(\RR)$.
 And in turn the set of real points in the open
$P$-orbit is obtained by throwing away the set of real points from
$P$-invariant but not $G$-invariant divisors (we denote them by $D_i$ and by $D_i^\circ$ the corresponding dense $P$-orbit). Any connected component
of $G$-orbit is defined by the the finite union of open
$P(\RR)$-orbits in the real topology. Without the loss of generality we
may assume that the real path $\gamma(t)$ joining any two points in the open
$G(\RR)^0$-orbit do not pass through the real points of the orbits of
codimension $\leq 2$ in particular we can assume that the path  lie in the open $P$-orbit for almost all points
and  intersect  $D_i^\circ$ in the finite set of points.  Assume that $P_\alpha$ raises $D_i^{\circ}$, then the neighborhood of the intersection point lies entirely
in $D_i^{\circ}$. In order to understand which $P(\RR)^\circ$-orbits form the connected component of the $G(\RR)$-orbit it is sufficient to study only the
case when the split semi-simple rank is $1$ and the closed or intermediate orbit have codimension $1$.

Let us fix a unique $P^-$-invariant point $y_0\in Y$, and let $N_{X/Y,y_0}$ be the fiber over $y_0$ of a normal bundle to $Y$. The open $\epsilon$-neighborhood which we denote by $N_{X/Y,\epsilon}$ of $Y$ in the real topology is diffeomorphic to the normal bundle $G*_{P^-} N_{X/Y,y_0}$ and the intersection  $N_{X/Y,\epsilon}\cap E_i$ defines the corresponding coordinate plane $H_i$ in $N_{X/Y,y_0}$. The connected components of
 $N_{X/Y,\epsilon}^\circ=N_{X/Y,\epsilon}\setminus \cup_i^{\rk X} E_i$ are marked by the components of $G(\RR)$-orbits to which they belong. Since the divisors $E_i$ form a complete intersection along $Y$ the element of $\pi_1(Y(\RR))$ defines the monodromy of the fiber $N_{X/Y,\epsilon}\setminus \cup_i^{\rk X} E_i$ over $y_0$ which is isomorphic to $ N_{X/Y,y_0}^\circ:= N_{X/Y,y_0}\setminus \cup_i^{\rk X} H_i $. Observe that the generators of $\pi_1(Y(\RR))$ are represented by the loops $\exp(\pi t Z_\alpha)y_0$ and $\exp(\pi t Z_\alpha)\in G(\RR)$. 
 For a point $x_0\in N_{X/Y,\epsilon}^\circ$ that defines a connected component we can represent the monodromy action by considering the path $\exp(\pi t Z_\alpha)x_0$ which lie entirely in the $G(\RR)^\circ$-orbit and by taking the corresponding connected component of $N_{X/Y,y_0}^\circ$. This actually proves that the monodromy action reduces to the action of $\exp(\pi Z_\alpha)$ on $N_{X/Y,y_0}^\circ$.


\begin{bibdiv}
  \begin{biblist}


    \bib{Bien}{article}{ author={Bien, Frédéric}, title={Orbits,
        multiplicities and differential operators}, conference={
        title={Representation theory of groups and algebras}, },
      book={ series={Contemp. Math.}, volume={145},
        publisher={Amer. Math. Soc., Providence, RI}, }, date={1993},
      pages={199--227}, doi={10.1090/conm/145/1216190},}

    \bib{BorelSerre}{article}{ author={Borel, Armand}, author={Serre,
        Jean-Pierre}, title={Théorèmes de finitude en cohomologie
        galoisienne}, journal={Comment. Math. Helv.}, volume={39},
      date={1964}, pages={111--164}, }

    \bib{BorelTits}{article}{ author={Borel, Armand}, author={Tits,
        Jacques}, title={Groupes réductifs}, journal={Inst. Hautes
        Études Sci. Publ. Math.}, number={27}, date={1965},
      pages={55--150}, }

    \bib{Brion}{article}{ author={Brion, Michel}, title={Quelques
        propriétés des espaces homogènes sphériques},
      journal={Manuscripta Math.}, volume={55}, date={1986},
      number={2}, pages={191--198}, }

    \bib{BBORB}{article}{ author={Brion, Michel}, title={On orbit
        closures of spherical subgroups in flag varieties},
      journal={Commentarii mathematici helvetici}, volume={76},
      date={2001}, number={2}, pages={263--299}, publisher={Springer},}

     \bib{TC}{article}{ author={  Cupit-Foutou, Stéphanie}  author={Timashev, Dmitry A.}
     title={Real orbits of complex spherical homogeneous spaces: the split case},
journal={arXiv:1909.04958}, }

    \bib{Kempf}{article}{ author={Kempf, George}, title={Instability
        in invariant theory}, journal={Ann. of Math. (2)},
      volume={108}, date={1978}, pages={299--316}, }

    \bib{Kimelfeld}{article}{ author={Kimelfeld, Boris},
      title={Homogeneous domains on flag manifolds},
      journal={J. Math. Anal. Appl.}, volume={121}, date={1987},
      number={2}, pages={506--588}, issn={0022-247X},
      review={\MR{872237}}, doi={10.1016/0022-247X(87)90258-7}, }
    
    \bib{inv.mot}{article}{ author={Knop, Friedrich}, title={The
        asymptotic behavior of invariant collective motion},
      journal={Inventiones mathematicae}, volume={116}, number={1},
      pages={309--328}, date={1994}, publisher={Springer}, }

    \bib{BORB}{article}{ author={Knop, Friedrich}, title={On the set
        of orbits for a Borel subgroup}, journal={Commentarii
        Mathematici Helvetici}, volume={70}, date={1995}, number={1},
      pages={285--309}, publisher={Springer}, }

    \bib{KK}{article}{ author={Knop, Friedrich}, author={Krötz,
        Bernhard}, title={Reductive group actions},
      journal={Preprint}, date={2016}, pages={62 pp.},
      arxiv={1604.01005}, }

    \bib{KSch}{article}{ author={Krötz, Bernhard},
      author={Schlichtkrull, Henrik}, title={Finite orbit
        decomposition of real flag manifolds},
      journal={J. Eur. Math. Soc. (JEMS)}, volume={18}, date={2016},
      number={6}, pages={1391--1403}, issn={1435-9855},}
  
    \bib{Matsuki}{article}{ author={Matsuki, Toshihiko}, title={Orbits
        on flag manifolds}, conference={ title={Proceedings of the
          International Congress of Mathematicians, Vol. I, II},
        address={Kyoto}, date={1990}, }, book={
        publisher={Math. Soc. Japan, Tokyo}, }, date={1991},
      pages={807--813}, }

    \bib{Popov}{article}{ author={Popov, V. L.}, title={Contractions
        of actions of reductive algebraic groups}, language={Russian},
      journal={Mat. Sb. (N.S.)}, volume={130(172)}, date={1986},
      number={3}, pages={310--334, 431}, issn={0368-8666},
      translation={ journal={Math. USSR-Sb.}, volume={58},
        date={1987}, number={2}, pages={311--335},
        issn={0025-5734},},}

    \bib{Rosenlicht}{article}{ author={Rosenlicht, Maxwell},
      title={Some basic theorems on algebraic groups},
      journal={Amer. J. Math.}, volume={78}, date={1956},
      pages={401--443}, }

    \bib{Vinberg}{article}{ author={Vinberg, Èrnest B.},
      title={Complexity of actions of reductive groups},
      journal={Funktsional. Anal. i Prilozhen.}, volume={20},
      date={1986}, number={1}, pages={1--13, 96}, }

    \bib{inv}{article}{ author={Popov, Vladimir L.} author={Vinberg,
        Èrnest B.}, title={Algebraic geometry IV}, pages={123--278},
      year={1994}, publisher={Springer}, }

    \bib{FUND}{article}{ author={Wiggerman, Mark}, title={The
        fundamental group of a real flag manifold},
      journal={Indagationes Mathematicae}, volume={9}, year={1998},
      number={1}, pages={141--153},}

  \end{biblist}
\end{bibdiv}
\end{document}